%% file: main_V2.tex
\documentclass[journal]{new-aiaa} %for conference papers
\usepackage[utf8]{inputenc}

\usepackage{graphicx}
\usepackage{amsmath}
\usepackage[version=4]{mhchem}
\usepackage{siunitx}
\usepackage{longtable,tabularx}
\usepackage{soul}
\input{PackagesCommands}
\usepackage[section]{placeins}

\usepackage[colorinlistoftodos,textwidth=1 in,textsize=footnotesize]{todonotes}

\title{Computational Modeling and Learning-Based Adaptive Control of Solid-Fuel Ramjets}
\author{
Gohar T. Khokhar\footnote{Postdoctoral Research Associate, Aerospace and Mechanical Engineering Department, Tucson, AZ 85721, AIAA Member.},
Kyle Hanquist\footnote{Assistant Professor, Department of Aerospace \& Mechanical Engineering, 1130 N. Mountain Avenue, University of Arizona, Tucson, AZ 85721, AIAA Associate Fellow.}
}
\affil{Department of Aerospace and Mechanical Engineering,
University of Arizona}
\author{
Parham Oveissi\footnote{Graduate Student, Department of Mechanical Engineering, 1000 Hilltop Circle, Baltimore, MD 21250.}, 
Alex Dorsey\footnote{Undergraduate Student, Department of Mechanical Engineering, 1000 Hilltop Circle, Baltimore, MD 21250.}, 
Ankit Goel\footnote{Assistant Professor, Department of Mechanical Engineering, 1000 Hilltop Circle, Baltimore, MD 21250.},
}
\affil{Department of Mechanical Engineering, University of Maryland, Baltimore County}
\begin{document}

\maketitle

% \listoftodos

\begin{abstract}

Solid-fuel ramjets offer a compact, energy-dense propulsion option for long-range, high-speed flight but pose significant challenges for thrust regulation due to strong nonlinearities, limited actuation authority, and complex multi-physics coupling between fuel regression, combustion, and compressible flow.
This paper presents a computational and control framework that combines a computational fluid dynamics model of an SFRJ with a learning-based adaptive control approach. 
A CFD model incorporating heat addition was developed to characterize thrust response, establish the operational envelope, and identify the onset of inlet unstart.
An adaptive proportional–integral controller, updated online using the retrospective cost adaptive control (RCAC) algorithm, was then applied to regulate thrust.
Closed-loop simulations demonstrate that the RCAC-based controller achieves accurate thrust regulation under both static and dynamic operating conditions, while remaining robust to variations in commands, hyperparameters, and inlet states.
The results highlight the suitability of RCAC for SFRJ control, where accurate reduced-order models are challenging to obtain, and underscore the potential of learning-based adaptive control to enable robust and reliable operation of SFRJs in future air-breathing propulsion applications.

% 
% Control of the combustion process under hypersonic conditions remains a challenging problem.  
% In this paper, we investigate the application of a data-driven, learning-based control technique to regulate a combustion process evolving inside a solid fuel ramjet to regulate the generated thrust under unknown operating conditions. 
% % 
% A computational model to simulate the combustion dynamics is developed by combining compressible flow theory with equilibrium chemistry. \todoUA{KH: I think we agreed to write the abstract last, right? Because this currently reflects the analytical theory + NASA CEA fidelity, which we could include here if needed. Should be ready to update now.} 
% % 
% The computational model is simulated to ascertain the combustion dynamics' stability and establish the engine's operational envelope.
% % 
% Based on retrospective cost optimization, an online learning controller is then integrated with the computational model to regulate the generated thrust. 
% Numerical simulation results are presented to demonstrate the robustness of the adaptive control system. 

\end{abstract}

\section*{Nomenclature}
{\renewcommand\arraystretch{1.0}
\noindent\begin{longtable*}{@{}l @{\quad=\quad} l@{}}
$E$ & Total energy per unit mass \\
$\bar{F}^\rmc$ & Convective flux vector \\
$\bar{F}^\rmv$ & Viscous flux vector \\
$M$ & Mach number \\
$\dot{m}$ & Mass flow rate \\
$T$ & Thermodynamic temperature \\
$p$ & Static pressure \\
$P$ & Thermodynamic pressure \\
$P_t$ & Total pressure \\
$P_r$ & Prandtl number \\
% $R_u$ & Universal gas constant \\
$\mathcal{R}$ & Residual vector (CFD) \\
$U$ & Conservative variables vector \\
$\tau$ & Thrust \\
$\bar{v}$ & Velocity vector \\
$\gamma$ & Specific heat ratio \\
$\kappa$ & Thermal conductivity \\
$\mu$ & Dynamic viscosity \\
$\rho$ & Density \\
$\bar{\bar{\tau}}$ & Viscous stress tensor \\
$I_n$ & $n \times n$ identity matrix \\
$\otimes$ & Kronecker product \\
$r_k$ & Commanded thrust \\
$y_k$ & Measured thrust \\
$z_k$ & Output error ($r_k - y_k$) \\
$\gamma_k$ & Accumulated output error \\
$u_k$ & Control input \\
$\theta_k$ & Controller gain vector \\
$K_{P,k}, K_{I,k}$ & Proportional and integral gains \\
$\Phi_k$ & Regressor matrix \\
$\overline w$ & Nominal heat flux \\
$w_k$ & Adaptive heat flux input \\
$K_w$ & Heat flux scaling factor \\
$P_0, N_1$ & RCAC hyperparameters \\
$R$ & Range (guidance law) \\
$\beta$ & Line-of-sight angle \\
$a_{z,c}$ & Commanded normal acceleration \\
$a_z$ & Measured normal acceleration \\
$\omega$ & Pitch rate \\
$\delta$ & Fin deflection angle \\
\end{longtable*}}

% \tableofcontents

\clearpage
\section{Introduction}

\iffalse 
what is the problem?
What prior work is done to solve the problem?
What is the current challenge?
What is done in this work to solve the challenge?
\fi 

% Ramjet engines, due to their stability and ability to generate high thrust for extended durations, are particularly well-suited for long-range operations at high speeds.  
% Because ramjets lack rotating turbomachinery, they are simpler to design, operate, and maintain compared to other air-breathing propulsion systems.
% % they are much smaller and simpler in construction than conventional jet engines.
% % 
% Ramjets can be classified based on the type of fuel they use: liquid-fuel ramjets (LFRJ) or solid-fuel ramjets (SFRJ).
% % 
% Due to the lack of turbopumps, fuel bladders, injectors, and associated plumbing, the SFRJ is much simpler than a comparably-sized LFRJ.
% Additionally, due to the higher volumetric energy density of solid fuel, the SFRJ has the potential for a greater range than a comparably-sized LFRJ.
% For an SFRJ, the combustion flame front spans the entire length of the fuel grain, which results in SFRJs being less likely than LFRJs to suffer from combustion instabilities.
% The solid fuel grain in an SFRJ offers the added advantage of being able to be stored as a solid circumferential module lined on the inner wall of the combustion chamber, eliminating the usual logistical concerns associated with liquid fuels.

%\textbf{Why study SFRJ?}
Ramjet engines are well-suited for sustained high-speed, long-range flight owing to their inherent stability and ability to generate continuous thrust over extended durations \cite{heiser1994hypersonic, sutton2011rocket}. 
The absence of rotating turbomachinery simplifies their design, operation, and maintenance relative to other air-breathing propulsion systems \cite{curran2001scramjet}.
Ramjets are typically categorized by fuel type as either liquid-fuel ramjets (LFRJs) or solid-fuel ramjets (SFRJs). The SFRJ architecture is significantly simpler than that of a comparably sized LFRJ, as it eliminates the need for turbopumps, fuel bladders, injectors, and associated plumbing \cite{gong2017numerical}. 
Furthermore, the higher volumetric energy density of solid propellants affords the SFRJ the potential for greater range than an equivalently scaled LFRJ \cite{gong2018combustion}.
In an SFRJ, the combustion flame front extends along the entire length of the solid grain, thereby reducing susceptibility to combustion instabilities \cite{li2020combustion, rashkovskiy2018combustion}. 
The solid grain can also be configured as a circumferential liner along the combustion chamber wall, enabling compact integration and storage while mitigating the logistical challenges typically associated with handling and transporting liquid fuels \cite{park2024trend}.

% \begin{figure}[h]
%     \centering
%     \includegraphics[width = 0.75\textwidth]{Figures/RAMJETSTATIONSGTK.png}
%     \caption{A typical SFRJ cross section with various stations.}
%     \label{fig:SFRJ_schematic}
% \end{figure} 

%\textbf{Challenges in SFRJ Operation?}
A solid-fuel ramjet (SFRJ) is an air-breathing propulsion system in which atmospheric oxygen serves as the oxidizer, while a hydrocarbon-based solid grain provides the fuel.
During operation, high-speed incoming air is compressed through the inlet, raising its pressure and temperature before entering the combustion chamber.
As the heated airflow passes over the exposed inner surface of the solid fuel grain, thermal feedback and mass transfer cause the fuel to regress, releasing gaseous pyrolysis products into the core flow.
These vapors mix with the compressed air and sustain a flame front that extends along the length of the grain. 
The combustion process produces high-temperature gases that are expanded through a convergent–divergent nozzle to generate thrust.
Unlike liquid-fueled systems, the fuel mass flow rate in an SFRJ is governed by the regression characteristics of the solid grain, which couple strongly with airflow conditions and combustion dynamics.

In an SFRJ, the high-speed, oxygen-rich intake flow reacts with the exposed fuel grain surface, releasing chemical energy that is converted into flow kinetic energy and ultimately into thrust. 
Stable operation requires that the thermodynamic state within the combustor—namely the pressure, temperature, and mass-flow rate—remain within a narrow band of conditions. 
If the airflow is insufficient, heat addition may not sustain the required thrust. 
In contrast, excessive airflow can lead to inlet unstart due to over-energization of the core flow or combustion blowoff, either of which can result in flame extinction and sudden thrust loss.
Predicting this stable operating envelope analytically is extremely difficult due to the complex, coupled multi-physics of solid-fuel combustion, turbulent mixing, and compressible reacting flow.
Reliable operation, therefore, demands regulation strategies that maintain the combustor state within acceptable limits and that remain robust to parametric uncertainty and external disturbances. 
From a control perspective, thrust regulation is further complicated by the underactuated and highly nonlinear dynamics of the SFRJ, where passive fuel regression and evolving chamber geometry preclude the use of conventional throttle mechanisms.

Thrust regulation of solid-fuel ramjets (SFRJs) has been an active area of research for more than four decades. 
Owing to the high cost and complexity of experimental testing, numerous high-fidelity computational tools have been developed to simulate the high-dimensional, multi-timescale, multi-physics flow environment inside an SFRJ \cite{netzer1977modeling,stevenson1981primitive,ben1999theoretical,Sun2009_simulation,wang2015numerical,schwer2018liquid,schwer2019progress}. 
These simulations have provided insight into the intricate flow interactions and have qualitatively characterized the steady-state response of SFRJs.
Early efforts to regulate thrust relied on algebraic relations between the air mass flow rate and the resulting thrust, derived from conservation of mass, momentum, and energy under quasi-static flow assumptions \cite{stevenson1981primitive,campbell1992solid,pelosi2003bypass,pei2013numerical,wang2015numerical,macleod2016review}. 
While computationally tractable, this approach requires precise knowledge of SFRJ physical parameters; measurement inaccuracies lead to erroneous inputs and, consequently, unreliable thrust predictions. 
More critically, because system dynamics are not captured in the algebraic framework, transient behavior cannot be predicted, and engine stability cannot be guaranteed.
To address these limitations, closed-loop control methods for thrust regulation began to emerge in the early 2000s \cite{Durali2005,durali2005velocity}. 
In particular, \cite{Durali2005} identified a linear dynamic model of the SFRJ and developed an adaptive controller to regulate thrust. 
However, controller performance was validated only against the simplified linearized model, rather than against a more comprehensive nonlinear representation of the SFRJ, limiting its practical applicability.

%\textbf{Why CFD?}
With advancements in modern computing capabilities, computational fluid dynamics (CFD) has been increasingly employed to investigate SFRJs \cite{Sun2009_simulation,wang2015numerical,nusca1990,kessler2022performance,goodwin2022simulating}.  
This progress has enabled the capture of additional SFRJ physics, including 3D effects, viscous effects, and chemical kinetics.
However, capturing all relevant phenomena in an SFRJ remains computationally expensive.  
Furthermore, while incorporating each of these physical effects enhances the predictive capabilities of the computational model and provides insights into the key interactions among various physical processes within the SFRJ, a higher-fidelity model is not necessarily beneficial for improving the control system. 
Additionally, since CFD models are structured as executable computer code, a high-fidelity computational model is typically unsuitable for controller design.
Theoretical and empirical tools for assessing a dynamic system's stability and transient characteristics in a loop with a control system are limited to a small class of nonlinear systems, often represented by ordinary differential equations.  
Nevertheless, high-fidelity CFD models are extremely valuable for data-driven, learning-based techniques to synthesize and stress test the control system.
This work thus investigates the application of an online, learning-based control design technique, called retrospective cost adaptive control (RCAC), to regulate the thrust generated by an SFRJ. 
RCAC has been recently demonstrated as a viable technique to synthesize an adaptive control system to regulate the thrust and prevent inlet unstart in a liquid-fuel scramjet engine \cite{goel2015scramjet,goel2018retrospective,goel2019output}.

% Retrospective cost adaptive control is based on retrospective cost optimization, where an auxiliary cost function, based on measured data and past inputs, is optimized to update the control law continually, and thus requires minimal modeling information \cite{santillo2010adaptive}. 
% In most applications, a simple first-order transfer function, which can be easily obtained by open-loop system simulation or from experimentally collected data, is sufficient.
% % 
% Within the context of RCAC, this transfer function is called the \textit{target model} since the RCAC algorithm updates the control law to drive an internal transfer function to the user-specified target model.
% % 
% It has been observed in numerous simulations and physical experiments that RCAC is insensitive to the choice of the target model. For example, the magnitude of the scalar that parameterizes the target model can be varied by as much as two orders of magnitude without adversely affecting the closed-loop performance. 

%\textbf{Why RCAC?}
Retrospective cost adaptive control is based on retrospective cost optimization, wherein an auxiliary cost function—constructed from measured data and past control inputs—is minimized to iteratively update the control law, thereby requiring minimal a priori modeling information \cite{santillo2010adaptive}. 
In most applications, a simple first-order transfer function, easily obtained from open-loop simulations or experimental data, is sufficient to characterize the system dynamics. 
Within the RCAC framework, this transfer function is referred to as the target model, as the algorithm adapts the controller to drive an internal transfer function toward the user-specified target. 
Several simulation and experimental studies have demonstrated that RCAC is robust to the choice of target model; for instance, the scalar parameter that defines the target model can vary by as much as two orders of magnitude without significantly degrading closed-loop performance.
This robustness is desirable for SFRJ applications, where accurate models are challenging to obtain due to the coupled multi-physics of solid-fuel combustion, turbulent mixing, and high-speed reacting flows.

% A key feature of RCAC that makes it particularly attractive to flow-control problems is its ability to adapt online. 
% Since RCAC requires only measured data (and not a system model) to optimize the controller, it can be directly integrated with a numerical simulation for controller hyperparameter tuning and stress testing.
% % 
% These features enable RCAC to be tuned using a computationally simple simulation of the system and then adapt appropriately to a more realistic simulation of the system or the physical system itself.
% Moreover, this work demonstrates the unique capability of RCAC to be tuned with a lower-fidelity model and the adapt to a higher-fidelity model without retuning. 
% % 
% Therefore, this work considers several fidelities to strike a balance, as shown in Figure \ref{fig:pyramid}, between resolving complex physics and computational cost, with the computational cost being limited as to the time needed to train and couple the controller.
% 

\begin{figure}[h]
    \centering
    \includegraphics[width = 0.495\textwidth]{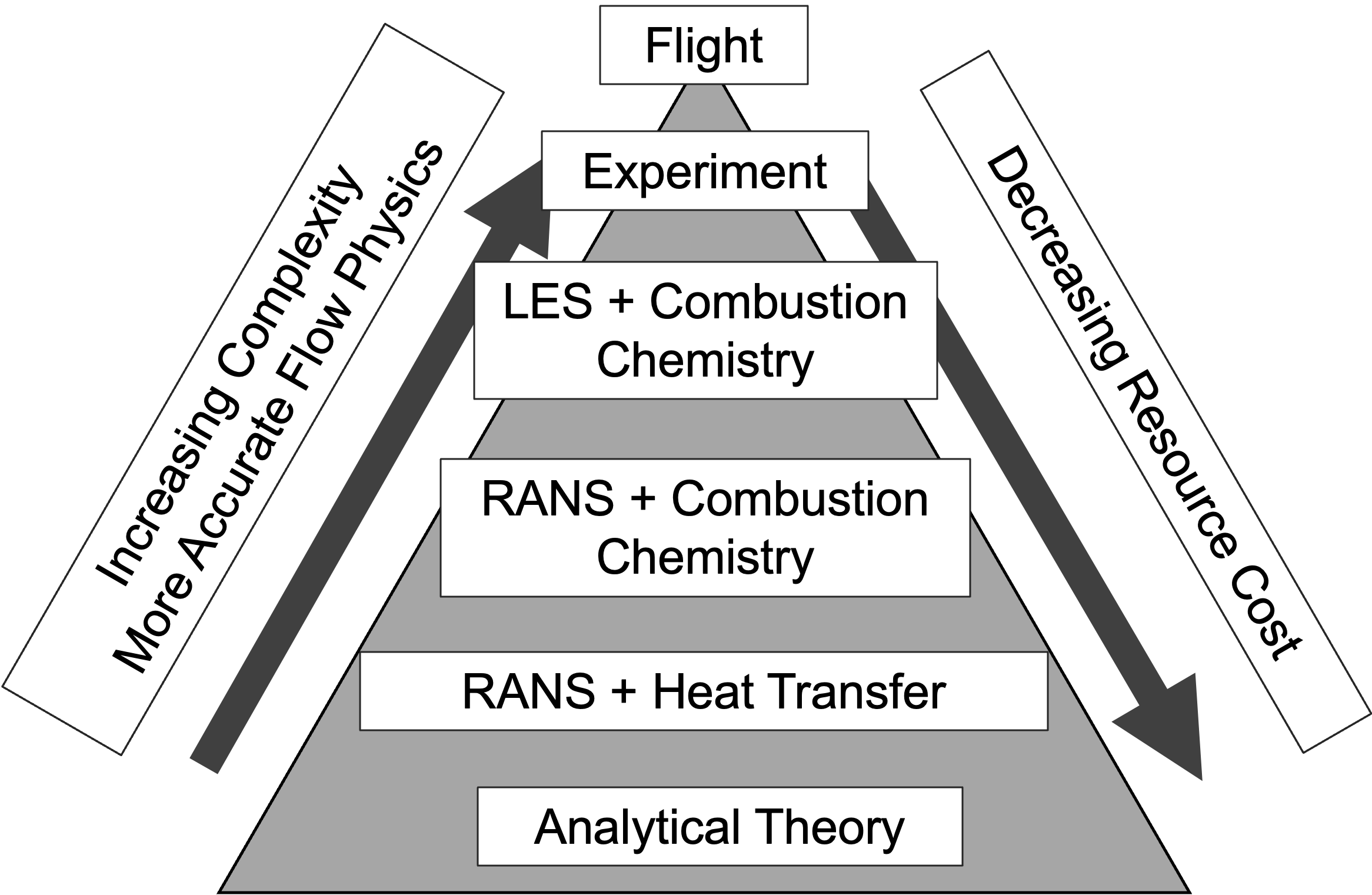}
    \caption{Different fidelities to investigate SFRJs where accuracy comes at the cost of resources.}
    \label{fig:pyramid}
\end{figure}

A key feature that makes RCAC particularly well-suited to flow-control problems is its ability to adapt online. 
Because RCAC optimizes the controller using only measured data, without reliance on an explicit system model, it can be directly integrated with numerical simulations for controller hyperparameter tuning and stress testing. 
This capability allows the controller to be tuned on a computationally inexpensive, low-fidelity simulation and then adapt appropriately when deployed on a higher-fidelity model or the physical system itself, without the need for retuning. 
Our ongoing efforts are focused on highlighting this unique capability of RCAC by considering multiple model fidelities, thereby striking a balance between capturing complex physics and managing computational cost, with the latter constrained primarily by the time required to train and couple the controller, as illustrated in Figure \ref{fig:pyramid}.

Our previous work \cite{oveissi2023learning,oveissi2025learning} used an analytical model based on equilibrium thermodynamics and chemistry to investigate the use of the data-driven learning technique for thrust regulation without reliance on the analytical model.
This approach, despite its efficiency and robustness, overlooks key SFRJ physics, such as geometry-specific and transient effects, which cannot be modeled using purely analytical theory due to the numerous empirical parameters involved.
This work focuses on constructing a computational approach that balances computational cost and fidelity to inform the control systems described later in Section \ref{sec:control}.
Preliminary results from this investigation were reported in \cite{oveissi2024adaptive,oveissi2025adaptive,khokhar2025investigation}.
% presents preliminary investigations into the application of retrospective cost adaptive control (RCAC) for regulating thrust in SFRJs.
% 
% This work focuses on constructing a computational approach that balances computational cost and fidelity to inform the control systems described in Section \ref{sec:control} and in detail in \cite{oveissi2023,oveissi2024}.
% Previous work utilized analytical theory to investigate the physical parameters of an SFRJ. This approach, despite its efficiency and robustness, overlooks key SFRJ physics, such as geometry specific and transient effects, which cannot be modeled using purely analytical theory due to the numerous empirical parameters involved.
% 
The main contributions of this work are:
\begin{enumerate}
\item the development of a CFD model with heat addition to capture combustion effects and its application to a realistic SFRJ geometry,
\item the design of an adaptive PI controller optimized using retrospective cost adaptive control (RCAC) without reliance on an explicit system model, and
\item the demonstration of robustness of the adaptive control system to variations in commands, hyperparameters, and operating conditions.
\end{enumerate}

This paper is organized as follows. 
Section \ref{sec:CFD} provides an overview of the governing equations used to model the flow and combustion in the SFRJ, describes the numerical solver and geometry considered in this work, and presents simulation results to establish the operational envelope of the SFRJ and to predict inlet unstart.
Section \ref{sec:control} describes the application of the retrospective cost adaptive control algorithm for thrust regulation in an SFRJ and presents closed-loop simulation results demonstrating successful performance under both static and dynamic operating conditions. 
The paper concludes with a discussion in Section \ref{sec:conclusions}.

% First, we provide preliminary investigations on a simplified SFRJ geometry using computational fluid dynamics (CFD) (Sec.~\ref{sec:CFD}).
% While this also does not include combustion chemistry, it is more dependent on the geometry of the SFRJ domain and not just the area ratios.
% Then we couple with the controller...\todoUMBC{Ankit/team: can you add a preview of your section here?}
% This is the first time...\todoUA{Brainstorm what is our one sentence contribution here.}

% \clearpage
\section{Modeling of SFRJ}
\label{sec:CFD}

\subsection{Computational Fluid Dynamics}

This work assumes that the flow is fully turbulent standard air, which is a suitable assumption for an SFRJ.
% \todo{AG: @Gohar, what does turbulent base flow mean? Why is it chosen to be investigated?}
To simulate the turbulent flow, the compressible Navier-Stokes equations are propagated, which in differential form are given by
\begin{equation}
     \mathcal{R}(U) = \frac{\partial U}{\partial t} + \nabla \cdot \bar{F}^{\rmc}(U) - \nabla \cdot \bar{F}^{\rmv}(U,\nabla U)   = 0
\end{equation}
where 
\begin{equation}
    U 
        \isdef 
            \matl 
                \rho \\ 
                \rho \bar{v} \\
                \rho E
            \matr 
            % \left \{  \rho, \rho \bar{v},  \rho E \right \}^\mathsf{T}
\end{equation}
is the conservative variable consisting of the fluid density $\rho,$
the velocity vector $\bar{v} \isdef \matl u & v & w \matr^\rmT,$ and 
the total energy per unit mass $E$.
% the conservative variables, $U$, are the working variables and given by:
% \begin{equation}
%     U = \left \{  \rho, \rho \bar{v},  \rho E \right \}^\mathsf{T}
% \end{equation}
% where $\rho$ is the fluid density, $\bar{v}$ is the velocity vector, and $E$ is the total energy per unit mass.
The convective flux $\bar{F}^{\rmc}$ and viscous flux $\bar{F}^{\rmv}$ are 
\begin{equation}
    \bar{F}^{\rmc}   
        \isdef
            % \left \{ \begin{array}{c} 
            \matl
            \rho \bar{v}  \\ \rho \bar{v} \otimes  \bar{v} + \bar{\bar{I}} p \\ \rho E \bar{v} + p \bar{v}
            \matr
            % \end{array} \right \}
    , \quad 
    \bar{F}^{\rmv} 
        \isdef
            % \left \{ \begin{array}{c} 
            \matl
                0 \\ \bar{\bar{\tau}} \\ \bar{\bar{\tau}} \cdot \bar{v} + \kappa \nabla T
                \matr
                % \end{array} \right  \}
\end{equation}
% 
%where $p$ is the static pressure, $\bar{\bar{\tau}}$ is the viscous stress tensor, $T$ is the temperature, $\kappa$ is the thermal conductivity, and $\mu$ is the viscosity.
%For the laminar flows of this work, the viscosity is simply the dynamic viscosity $\mu_\rmd$, which is assumed to satisfy Sutherland's law as a function of temperature.
% 
% where $\rho$ is the fluid density, $\bar{v}=\left\lbrace u, v, w \right\rbrace^\mathsf{T}$ is the flow speed in Cartesian system of reference, $E$ is the total energy per unit mass, 
where $p$ is the static pressure, $\bar{\bar{\tau}}$ is the viscous stress tensor, $T$ is the temperature, $\kappa$ is the thermal conductivity, and $\mu$ is the viscosity, which is assumed to satisfy Sutherland's law as a function of temperature. 
Note that the viscous stress tensor $\bar{\bar{\tau}}$ can be written as 
\begin{align}
    \bar{\bar{\tau}}= \mu \left ( \nabla \bar{v} + \nabla \bar{v}^{T} \right ) - \mu \frac{2}{3} \bar{\bar I} \left ( \nabla \cdot \bar{v} \right ).
\end{align}
where the operator $\nabla$ is the gradient vector, 
and $\bar{\bar I}$ is $3 \times 3$ identity matrix. 
% \todo{AG: @Gohar, cite a reference for the equation above}

Assuming a perfect gas with a ratio of specific heats $\gamma$ and specific gas constant $R$, the system of equations is closed by using
\begin{align}
p = (\gamma-1) \rho \left [ E - 0.5(\bar{v} \cdot \bar{v} ) \right ].    
\end{align}
In turbulent flows, to solve the Reynolds-averaged Navier-Stokes (RANS) equations, we use the Boussinesq hypothesis, which states that the effect of turbulence can be represented as an increased viscosity. 
In this work, all flow simulations are assumed to be fully turbulent and are modeled using the $k - \omega$ SST turbulence model.

%\textbf{Numerical solver.}
The CFD software used in this work is SU2, a computational analysis and design package developed for solving multiphysics analysis and optimization problems on unstructured mesh topologies \cite{economon2016a}.
SU2 employs a median-dual finite-volume method to discretize the governing equations.
Convergence is assessed by monitoring the root-mean-square residuals of mass and energy at cell centroids, as well as the global mass imbalance between the inlet and outlet, ensuring conservation of mass and energy throughout the domain.
For the present study, convective fluxes are discretized using the Jameson–Schmidt–Turkel (JST) scheme \cite{liou1993}. 
Additional details regarding the governing equations and numerical methods implemented in SU2 are provided in \cite{economon2016a}.

Combustion is modeled as simplified heat addition originating from the wall or fuel-grain surface.
% as illustrated in Figure \ref{fig:SFRJ_schematic}. 
Although this approach is less sophisticated than resolving detailed chemical kinetics, it captures the primary effect of combustion by representing the enthalpy change associated with the heat of reaction.
The heat addition is imposed along an isobaric path, thereby directly influencing the flow enthalpy. 
Furthermore, the finite-volume formulation employed in the simulations is inherently conservative, ensuring the conservation of mass, momentum, and energy fluxes throughout the computational domain.

\begin{figure}[h]
    \centering
    \includegraphics[width=0.9\columnwidth]{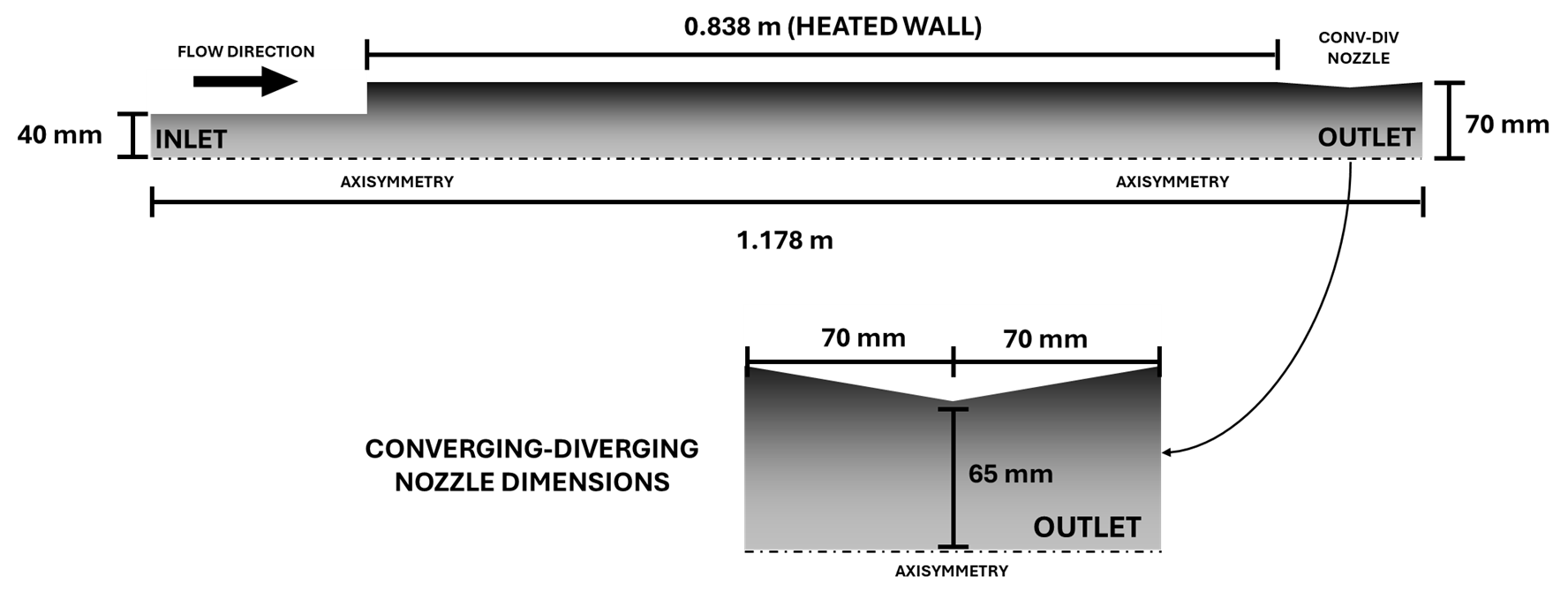}
    \caption{Simplified SFRJ geometry}
    \label{fig:Truncated_vs_Full}
\end{figure}

The present work considers a simplified backward-facing step geometry inspired by the full SFRJ configuration studied in Ref. \cite{li2021numerical}. 
Such fundamental geometries are widely employed in SFRJ simulation studies \cite{bojko2024investigating, deboskey2024analysis}. 
The geometry used in this investigation, illustrated in Figure \ref{fig:Truncated_vs_Full}, represents a truncated version of the full SFRJ geometry from Ref. \cite{li2021numerical}.
While retaining the original dimensions used in Ref. \cite{li2021numerical}, it comprises an inlet channel, a combustor, and an exit converging-diverging nozzle.
% \todo{@Gohar, original dimensions of what?}
The inlet channel has a diameter of $80 \ \rm mm$ and a length of $0.2 \ \rm m$. 
The combustor is $ 140 \ \rm mm$ in diameter and $0.838 \ \rm m$ in length.
The exit converging–diverging nozzle is symmetric, with a length of $140 \ \rm mm$ and a throat diameter of $130 \ \rm mm$. 
The total length of the SFRJ geometry is $1.178 \ \rm m$.
To model combustion in the SFRJ, a heat flux is imposed along the heated wall section shown in Figure \ref{fig:Truncated_vs_Full}, which has a length of $0.838 \ \rm m$.

% \todo{@Gohar, please change the text in the figure to "Simplified SFRJ geometry." Note the lowercase. DONE}
% \todo{Please mark inlet plane and output plane in the first figure. DONE}

\begin{figure}[h]
    \centering
    \includegraphics[width=\columnwidth]{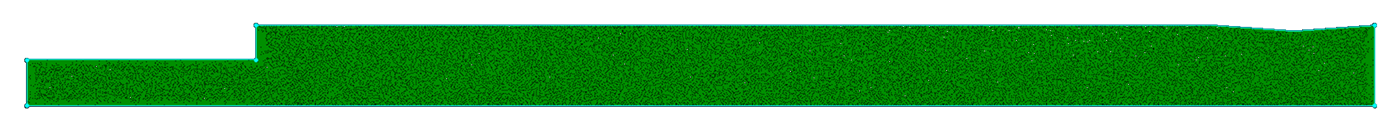}
    \caption{Computational mesh consisting of approximately 60,000 unstructured cells.}
    \label{fig:SFRJmesh}
\end{figure}

To reduce computational cost, the axisymmetric nature of the SFRJ configuration is exploited by simulating only half of the two-dimensional cross-section through the center of the full three-dimensional geometry. This is performed by simulating the domain as 2D axisymmetric where one dimension is the axial distance ($X$ is later figures) and the other is the distance from the centerline ($Y$ in later figures).
The computational domain, illustrated in Figure \ref{fig:SFRJmesh}, comprises approximately $60,000$ unstructured cells. A grid convergence study for this mesh is provided in Sec.~\ref{sec:grid}.

\begin{figure}[h]
    \centering
    \includegraphics[width=0.9\textwidth]{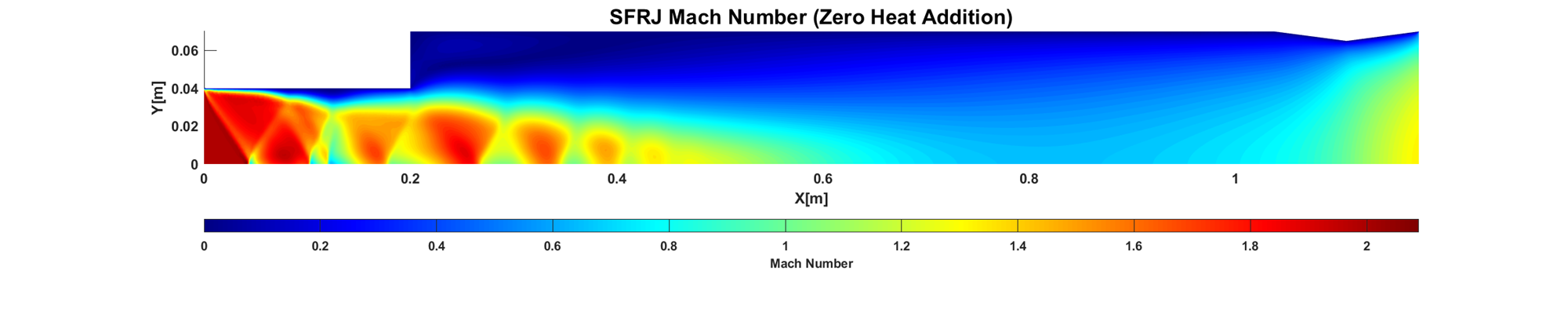}
    \caption{Mach number contour for SFRJ with no heat addition}
    \label{fig:MachContourPlot}
\end{figure}

At the inlet of the SFRJ, a velocity of 695 m/s, a static pressure of 100,000 Pa, and a static temperature of 300 K are prescribed. 
These conditions were chosen to achieve an inlet Mach number of 2, while the exit was prescribed as supersonic.
% \todo{Gohar, do you mean that the V, P, and T combination is equivalent to Mach 2? ADDRESSED}
% \todo{Dont we need to specify exit conditions as well to compute the steady state? ADDRESSED}
Figure \ref{fig:MachContourPlot} shows the constant Mach contours without any head addition at steady state. 
Note that the flow is supersonic at the inlet, subsonic in the combustor, and supersonic at the exhaust.
% These observations are consistent with the physical observations from a solid fuel ramjet operating in similar conditions. \todo{Need citation to support this.}

% Baseline Mach number contours, with no heat addition, for the  SFRJ are shown in Figure \ref{fig:MachContourPlot}. 

% Note the supersonic inlet, the subsonic combustor region, and the supersonic exhaust.
% All of the flow conditions are consistent with the expected steady state operating conditions of a solid-fuel ramjet. 

\begin{figure}[h]
    \centering
    \includegraphics[width=3.25 in]{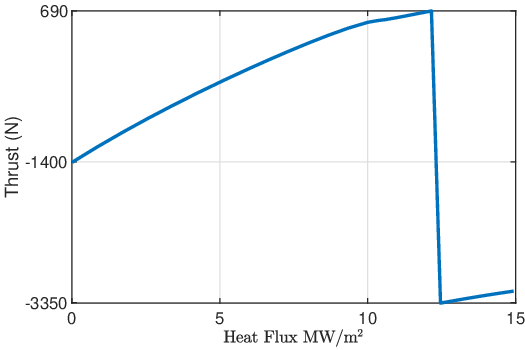}
    \caption{Thrust at various Heat Flux Values}
    \label{fig:open_loop_thrust_heatflux_0_1_8}
\end{figure}

Next, we consider the case of the heated wall, which is shown in Figure \ref{fig:Truncated_vs_Full}.
In particular, a constant heat flux is added to simulate a heated wall. 
In this work, we set the heat flux to be $\{2, 4, 6, 8, 10, 12, 14, 16 \} \times 10^6$ $\rm Watts/m^2$ and let the SFRJ reach a steady state.
% The resulting thrust increase is documented for each incremental rise in heat flux. 
Using mass-averaged quantities at the inlet and outlet, the thrust, $\tau$, generated by the SFRJ is calculated by applying the principle of momentum conservation to the control volume, that is, 
\begin{align}
\tau = \dot{m} (v_{\text{outlet}} - v_{\text{inlet}}) + (P_{\text{outlet}} - P_{\text{inlet}}) A_{\text{outlet}}.
\label{eq:thrust}
\end{align}
% Next, we compute the relative thrust defined as the difference between the thrust generated with a constant heat flux and no heat flux.
Figure \ref{fig:open_loop_thrust_heatflux_0_1_8} shows the thrust generated by the SFRJ as the heat flux is increased. 
Note the sudden loss of thrust beyond $12 \times 10^6$ $\rm Watts/m^2$ due to what we are describing as unstart, which is discussed more in the following section.
% steady-state change in thrust (delta thrust) generated by the SFRJ for increasing wall heat flux values.
% Note the nonlinear input-output relation between the delta thrust and heat flux, and the sudden loss of thrust around \(12 \, \text{million W/m}^2\).

% \begin{figure}[h]
%     \centering
%     \includegraphics[width=3.25 in]{Figures/QvsT_CFD.png}
%     \caption{Steady-state thrust generated by the SFRJ for heat flux values.}
%     \label{fig:open_loop_thrust_heatflux_0_1_8}
% \end{figure}
% \todoUMBC{@Alex, lets show the absolute thrust here, not the delta thrust. Remove legend}

% \todoUMBC{Alex, we need an 8x2 figure here showing Mach in the first column and the pressure in the second column for all values of heat flux that you used.}

\begin{figure}[h]
    \centering
    \includegraphics[width=3.25 in]{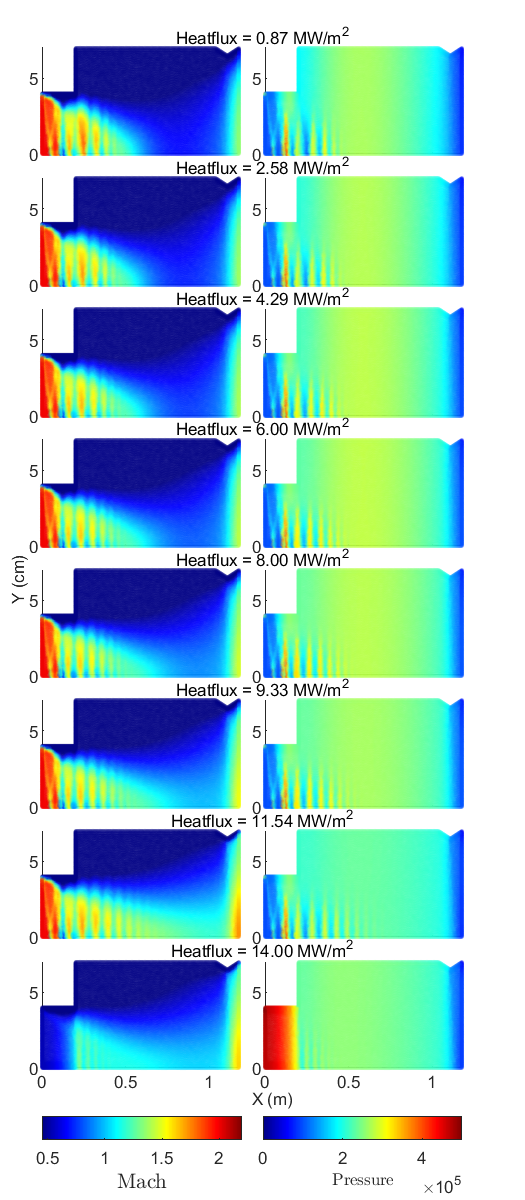}
    \caption{Mach and Pressure contours at various heat flux values}
    \label{fig:MachAndPressure8by2}
\end{figure}

% \todo{@Gohar, edit the folllowing discussion about mesh sensitivity DONE}

\subsection{Unstart}

% The sudden drop in thrust is due to thermal choking in the combustor region, leading to engine unstart. In a subsonic flow, heat addition increases the Mach number, causing the flow in the combustor to approach sonic conditions. Once the flow becomes choked, further heat addition cannot increase kinetic energy; instead, it raises the internal energy, manifesting as increased temperature and pressure. 
% This rise in temperature and pressure reduces the pressure differential between the inlet and outlet of the SFRJ, disrupting the supersonic condition at the inlet, thus resulting in a dramatic reduction in thrust. 
% Figure \ref{fig:SFRJMachstacked} shows the Mach contours before and after engine unstart. 
% Note the sonic conditions now present in the combustor region when compared to Figure \ref{fig:MachContourPlot}. 
% Figure \ref{fig:SFRJPressurestacked} shows the static pressure contours before and after engine unstart.
% Note the sudden change in the inlet pressure.
% The inlet pressure increase is directly responsible for the loss in computed thrust from \eqref{eq:thrust}.
The sudden drop in thrust is caused by thermal choking in the combustor region, resulting in an engine unstart. 
In subsonic flow, heat addition increases the Mach number, causing the flow within the combustor to approach sonic conditions.
Once choking occurs, further heat addition cannot increase the kinetic energy of the flow; instead, it raises the internal energy, resulting in higher temperature and pressure.

\begin{figure}[h]
    \centering
    % Subfigure (a)
    \begin{subfigure}[b]{0.49\textwidth}
        \centering
        \includegraphics[width=\textwidth]{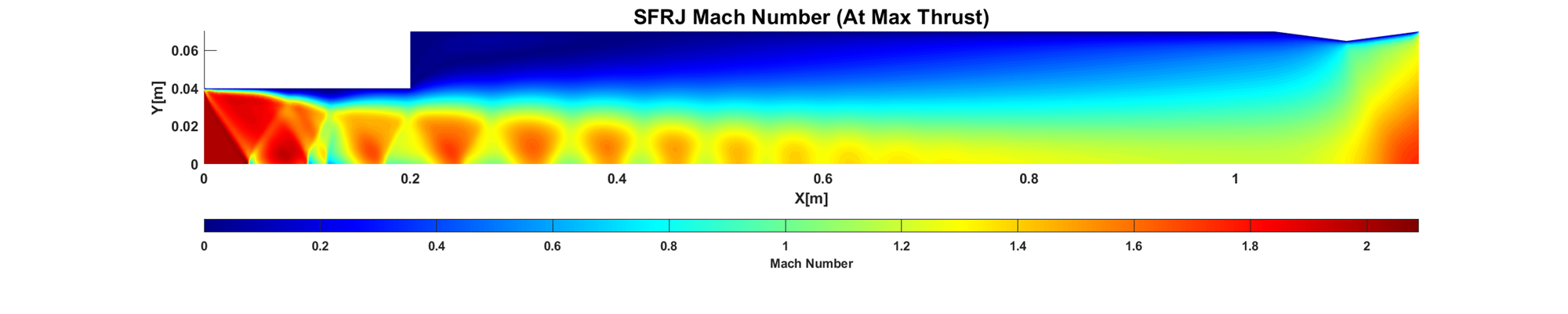} % Replace with your image file
        \caption{Mach contour before unstart.}
        \label{fig:SFRJMaxThrust}
    \end{subfigure}
    \\
    \begin{subfigure}[b]{0.49\textwidth}
        \centering
        \includegraphics[width=1\textwidth]{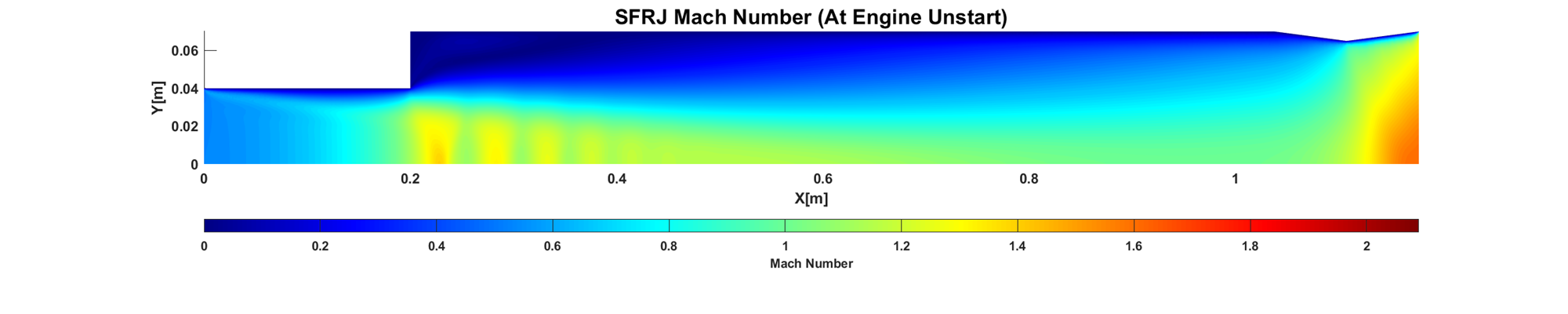} % Replace with your image file
        \caption{Mach contour after unstart.}
        \label{fig:SFRJUnstart}
    \end{subfigure}
    \\
    \begin{subfigure}[b]{0.49\textwidth}
        \centering
        \includegraphics[width=\textwidth]{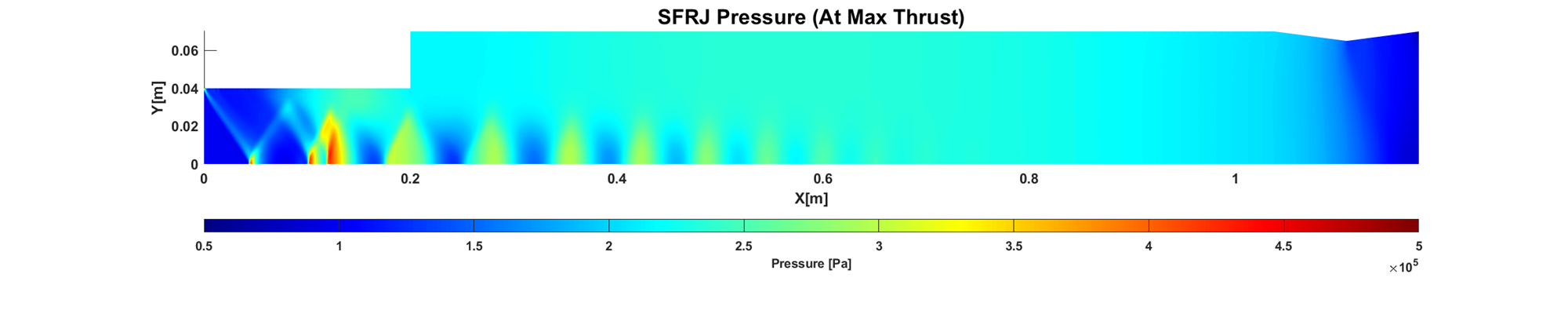} % Replace with your image file
        \caption{Static pressure contour before unstart.}
        \label{fig:SFRJPressMaxThrust}
    \end{subfigure}
    \\
    \begin{subfigure}[b]{0.49\textwidth}
        \centering
        \includegraphics[width=1\textwidth]{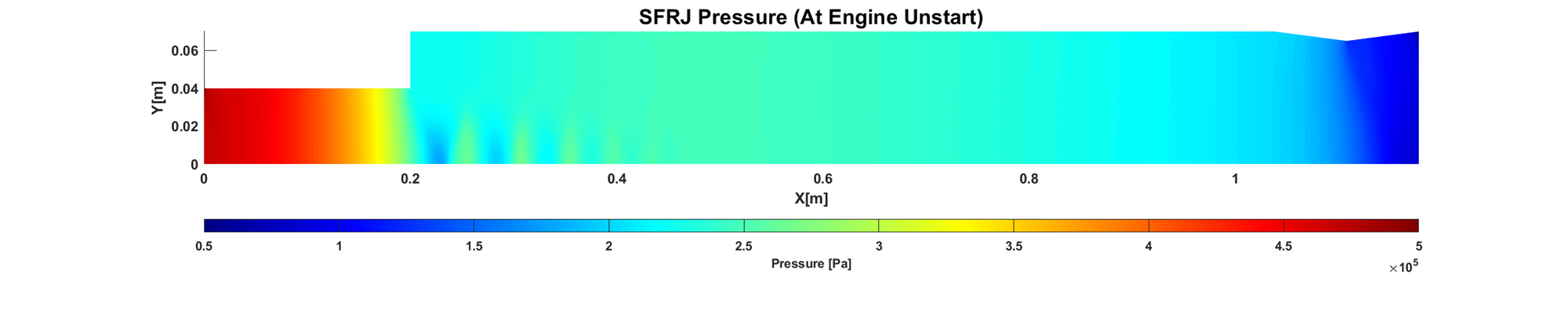} % Replace with your image file
        \caption{Static pressure contour after unstart.}
        \label{fig:SFRJPressUnstart}
    \end{subfigure}

    % Overall caption
    \caption{Thermodynamic contours before and after engine unstart in the SFRJ.}
    \label{fig:SFRJMachstacked}
\end{figure}

This rise in temperature and pressure diminishes the pressure differential between the inlet and outlet of the SFRJ, disrupting the supersonic condition at the inlet and causing a significant reduction in thrust. 
Figures \ref{fig:SFRJMaxThrust} and \ref{fig:SFRJUnstart} show the Mach number contours before and after engine unstart, highlighting the emergence of sonic conditions in the combustor compared to Figure \ref{fig:MachContourPlot}.
Figures \ref{fig:SFRJPressMaxThrust} and \ref{fig:SFRJPressUnstart} show the static pressure contours before and after unstart, revealing a marked increase in inlet pressure. 
This rise in inlet pressure directly contributes to the reduction in computed thrust, defined by \eqref{eq:thrust}.

\begin{figure}[h]
    \captionsetup{justification=centering}
    \centering
    \begin{subfigure}{0.32\textwidth}
        \centering
        \includegraphics[width=\linewidth]{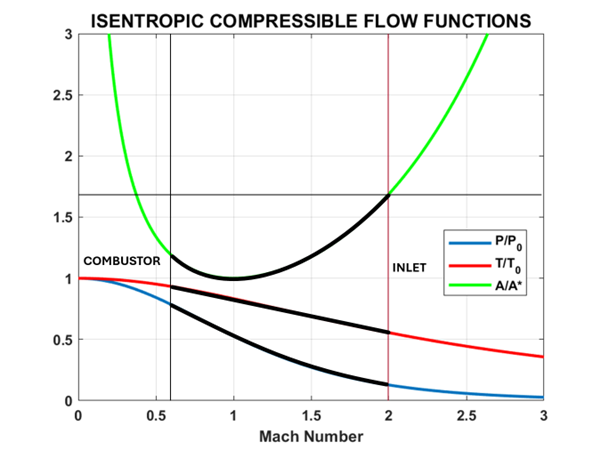}
        % \caption{Inlet and combustor states at baseline thrust (before heat addition)}
        \caption{Without heat addition.}
        \label{fig:a}
    \end{subfigure}
    % \hfill
    \begin{subfigure}{0.32\textwidth}
        \centering
        \includegraphics[width=\linewidth]{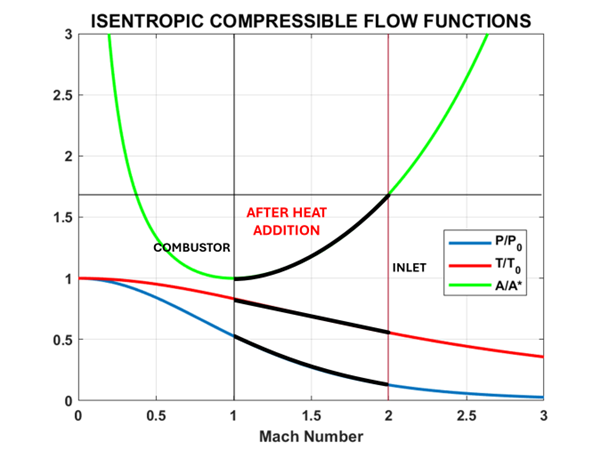}
        % \caption{Inlet and combustor states at max thrust (right before engine unstart)}
        \caption{Before engine unstart.}
        \label{fig:b}
    \end{subfigure}
    % \hfill
    \begin{subfigure}{0.32\textwidth}
        \centering
        \includegraphics[width=\linewidth]{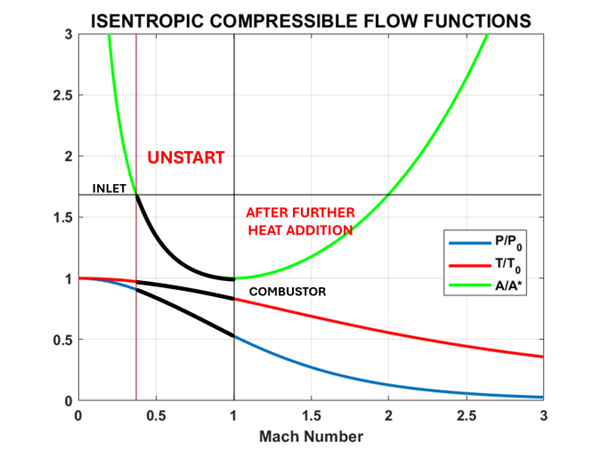}
        \caption{After engine unstart.}
        \label{fig:c}
    \end{subfigure}
    \caption{Approximation of SFRJ inlet and combustor states using compressible flow functions.}
    \label{fig:IsentropicMain}
\end{figure}

The engine unstart prediction from the CFD simulations can be corroborated with analytical theory by documenting inlet and combustor states using isentropic compressible flow functions.
For the current simulations, the flow is modeled as viscous; however, the boundary layers are small enough to obtain reasonable approximations of thermodynamic states using isentropic flow functions. 
Figure \ref{fig:IsentropicMain} plots the evolution of the inlet and combustor thermodynamic states as heat is added. 
The static-to-total pressure, static-to-total temperature, and area-to-sonic area ratios are shown.
Figure \ref{fig:a} depicts the baseline thrust scenario before any heat addition occurs.
Note the subsonic combustor state and supersonic inlet state. 
Figure \ref{fig:b} shows the evolution of the combustor state as the SFRJ reaches its max thrust output. 
At this point, the combustor reaches the sonic condition, that is, $M =1$.
Adding any extra heat to the combustor beyond this does not change the Mach number of the combustor, since heat addition to a compressible flow, subsonic or supersonic, will only result in the flow tending towards the sonic condition.
Extra heat addition beyond the state depicted in Figure \ref{fig:b}, results in the internal energy of the system directly increasing, resulting in a temperature rise. 
Since the combustor thermodynamic state, including temperature, becomes immutable at the sonic condition, the only possible temperature rise the system can accommodate corresponds to the alternate inlet subsonic area to sonic area ratio.
The sudden shift in the inlet condition from sonic to subsonic allows for the necessary temperature rise to occur while maintaining the inlet's cross-sectional area. 
From Figure \ref{fig:c}, the sudden new inlet temperature needed to satisfy energy conservation occurs at a higher static pressure compared to Figure \ref{fig:b}. 
This sudden transition to a higher static pressure in Figure \ref{fig:c} is ultimately responsible for the sudden diminished thrust the SFRJ experiences at engine unstart.

% The computational domain shown in Figure \ref{fig:SFRJmesh} consists of approximately 60,000 unstructured cells. A grid convergence study, shown in Figure \ref{fig:GridDep}, was conducted using meshes with 30k, 60k, 90k, and 120k cells to determine an appropriate resolution for simulating the SFRJ. The study showed that the predicted thrust resulting from heat addition converged at a cell count of 60k. However, the predicted location of the SFRJ unstart was found to be sensitive to grid resolution, with higher cell counts predicting a delayed onset of unstart. Since the controller operates only in the pre-unstart regime, it was decided that the 60k cell mesh would be sufficient for this study. 

%pre-unstart heat fluxes are grid converged, and computational cost considerations (e.g., CFD is coupled to the controller), 
%, further investigation into the grid dependency of the unstart location was not pursued. 
%\todo{@Gohar, can you add more detail to the mesh details? Explain why 60,000 is appropriate and what scheme is used to generate this. Just write as much as you can, we will refine it as we iterate over the paper; KH: yes, I agree we need to provide more details and why this mesh is suitable. We could quickly generate meshes with different qualities and show a mesh vs. thrust convergence study. We need to show sort of justification for the mesh as I've had papers rejected that do not include this. DONE}

% More details about the domain geometry can be found in Sec.~\ref{sec:cfdresult}.

\subsubsection{Grid convergence study}
\label{sec:grid}
\begin{figure}[h]
    \centering
    \includegraphics[width=0.5\linewidth]{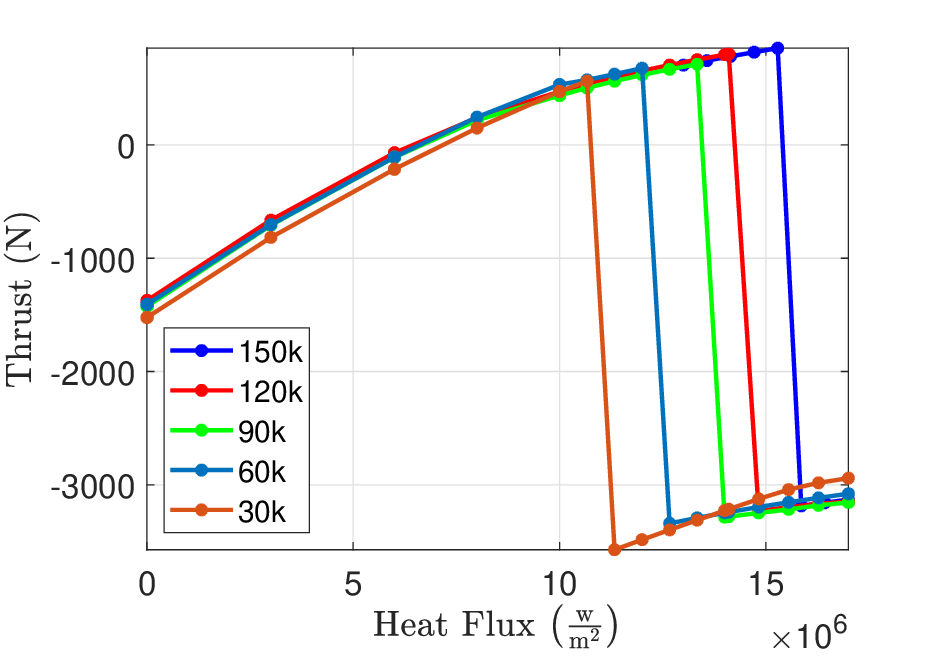}
    \caption{Grid dependence of SFRJ predicted thrust due to heat addition}
    \label{fig:GridDep}
\end{figure}

Finally, a grid convergence study was performed using meshes with $30,000$, $60,000$, $90,000$, and $120,000$ cells, as shown in Figure \ref{fig:GridDep}, to identify the sufficient mesh resolution to simulate the SFRJ. 
The results indicated that the predicted thrust due to heat addition converged at a mesh size of $60,000$ cells. 
However, the predicted location of SFRJ unstart exhibited sensitivity to grid resolution, with finer meshes predicting a delayed onset of unstart. 
Since the controller operates exclusively in the pre-unstart regime, the $60,000$-cell mesh was deemed sufficient for this study.

\section{Learning-based Control System}
\label{sec:control}
The previous section developed a predictive model to compute the thrust generated by the SFRJ and to identify the unstart region to be avoided. 
This section addresses thrust regulation using an adaptive PI controller that is continually updated by the retrospective cost adaptive control (RCAC) algorithm. 
The adaptive controller is then applied to regulate thrust across a range of operating conditions, and its robustness is demonstrated through numerical simulations.

% Now that we have a computational model of the SFRJ that predicts thrust based on inputs (\todoUMBC{Add inputs used here.}) and has a region of unstart to avoid, the goal now is to control it.\todoUMBC{Does this paragraph sound okay? Trying to transition between the topics.} 

% \begin{enumerate}
%     \item Static tests
%     \begin{todolist}
%         \item [\done] Step command
%         \item [\done] Doublet
%         \item [\done] Effect of hyperparameters
%     \end{todolist}

%     \item Dynamic tests
%     \begin{todolist}
%         \item Step command
%         \item Doublet
%         \item Effect of hyperparameters
%     \end{todolist}

%     % \item Control near unstart

%     % \item Neural network map of thrust with pressure measurements
% \end{enumerate}

\subsection{Algorithm}
\begin{figure}[h]
    \centering
    \resizebox{3.25 in}{!}
    {
    \begin{tikzpicture}[auto, node distance=2cm,>=latex',text centered]
    
        \node at (-3,0) (reference) {$r_k$};
        \node[sum, right = 3 em of reference] (sum) {};
        \draw[->] (sum.east) -- +(0.5,0) -- +(0.5,-1) -- +(1.25,-1)
                    -- +(3,1);
        \node [smallblock, fill=green!20, right = 3 em of sum] (Controller) 
        {$\begin{array}{c}{\rm Adaptive \ PI} \\ {\rm Controller}\end{array}$};

        \node [smallblock, fill=green!20, right = 3 em of Controller, minimum height = 1.2cm] (Nonlinear) {Map};
        % {$\begin{array}{c}{\rm Nonlinear} \\ {\rm Map}\end{array}$};

        % \node [smallblock, fill=green!20, right = 3 em of Nonlinear] (ZOH) {ZOH};

        \node [smallblock, right = 3 em of Nonlinear] (Plant) 
        % {SFRJ};
        {$\begin{array}{c}{\rm Computational} \\ {\rm SFRJ \ Model}\end{array}$};
        
        % \node [smallblock,  fill=green!20, right = 3 em of Plant] (A2D) {A2D};
        
        % , minimum height=3em, text width=1.6cm

        % \node at {(Plant)+(3,0)}  (output) {$y$};
        % \node[right = 5 em of Plant] (output) {$y$};
        
        \draw[->] (reference) -- node[xshift = 1.5em, yshift = -1.7em]{$-$} (sum);
        \draw[->] (sum) -- node[xshift = 0em, yshift = .2em]{$z_k$} (Controller);
        \draw[->] (Controller)-- node[xshift = 0em, yshift = .2em]{$u_k$} (Nonlinear);
        \draw[->] (Nonlinear)-- node[xshift = 0em, yshift = .2em]{$w_k$} (Plant);
        % \draw[->] (ZOH) -- node[xshift = 0em, yshift = .2em]{$\overline w$} (Plant);
        % \draw[->] (Plant)-- node[xshift = 0em, yshift = .2em]{$\overline{y}$} (A2D);
        \draw[->] (Plant.east) -- node[xshift = 0em, yshift = .2em]{$y_k$} +(1,0) --  +(+1,-1.5)
                    -| (sum.south);
        % \draw[->] (Plant.north)+(0,1) node[xshift = 0.5em, yshift = -0.5em]{$d$} -- (Plant.north);

    \end{tikzpicture}
    }
    % \vspace{-2em}
    \caption{Control architecture to regulate the thrust generated by the SFRJ. }
    \label{fig:BSL}
\end{figure}
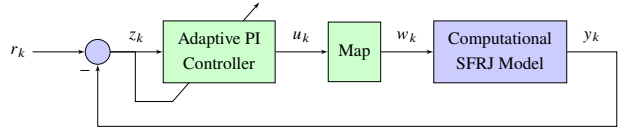

This section provides a brief overview of the learning-based control system used to regulate the thrust generated by the SFRJ model introduced in the previous section.
Figure \ref{fig:BSL} illustrates the closed-loop feedback architecture.
The control system comprises an adaptive proportional–integral (PI) controller, whose gains are continuously updated using the retrospective cost adaptive control (RCAC) algorithm.
In addition, an affine function is used to translate the RCAC-generated control signal into the corresponding heat flux input.
The mapping is designed to ensure that the magnitude of the adaptive control signal remains on the order of $\SO(1)$, thereby preserving the numerical stability of the optimization routine within the RCAC framework.
% The map is chosen to ensure that the optimization problem solved internally in RCAC remains well-conditioned. 
% The nonlinearity is chosen based on the open-loop simulations of the SFRJ model that showed the nonlinear relationship between the steady-state heat flux and the steady-state thrust output of the SFRJ. 

The adaptive PI controller can be written as
\begin{align}
    u_k = K_{\rmP,k} z_k + K_{\rmI,k} \gamma_k,
    \label{eq:PI_control}
\end{align}
where $z_k$ is the output error defined as the difference between the commanded thrust $r_k$ and the measured thrust output $y_k,$ that is, $z_k \isdef r_k - y_k,$
$\gamma_k$ is the accumulated output error given by
\begin{align}
    \gamma_k
        \isdef 
            \sum_{i=0}^k z_i, 
\end{align}
and the scalars $K_{\rmP,k} $ and $K_{\rmI,k}$ are the proportional and integral gains optimized by the RCAC algorithm at step $k. $
Note that the integral signal $\gamma_k$ can be computed recursively as
$
    \gamma_{k+1} = \gamma_k + z_{k+1}.
$

Next, the adaptive PI control law \eqref{eq:PI_control} is reformulated in the regressor form as
\begin{align}
    u_k = \Phi_k \theta_k, 
    \label{eq:PI_control_ref}
\end{align}
where 
\begin{align}
    \Phi_k 
        \isdef 
            \matl 
                z_k & \gamma_k
            \matr, 
    \quad
    \theta_k 
        \isdef 
            \matl 
                K_{\rmP,k} \\
                K_{\rmI,k}
            \matr,
\end{align}
where the regressor matrix $\Phi_k$ contains the measured data and the controller gain vector $\theta_k$ is optimized by the RCAC algorithm described in \cite{oveissi2023learning,poudel2023learning}.

Finally, the heat flux $w_k$ is given by 
\begin{align}
    w_k
        =
            \overline w  + K_w {u_k},
    \label{eq:nonlinear_map}
\end{align}
% where $\overline w$ is the nominal heat flux and 
% $K_w$ is the scaling factor that is chosen based on the open-loop simulations as discussed below. 
% Note that the units of the heat flux are $\rm W/m^2.$
% % The choice of $k_w$ will be discussed in Section \ref{xyz}.
% In this work, the nominal heat flux $\overline w = 10 \times 10^6$ $\rm W/m^2.$ 
% 
% \todoUMBC{@Parham, add values of wbar and Kw and write a short description of how you got them.}
% \todoUMBC{@Parham, finalize the nonlinear map description.}
where the nominal heat flux $\overline{w}$ and the scaling factor $K_w$ are set to $10 \times 10^6$ $ \rm{W/m^2}$ and $10^6$, respectively.
These values are selected based on open-loop simulations. 
Specifically, $\overline{w}$ is chosen as a baseline heat flux that lies approximately midway between the value that yields zero net thrust and the value that produces maximum net thrust before engine unstart.
The scaling factor $K_w$ is selected to ensure that the adaptive control signal $u_k$ remains on the order of unity, that is, $\mathcal{O}(1)$, throughout the simulation. 
This choice promotes the numerical stability and reliability of the optimization algorithm employed in the RCAC method.

\subsection{Static Operating Conditions}

In this section, we evaluate the performance of the adaptive controller under static operating conditions.
These tests are inspired by benchtop experiments typically conducted in a laboratory setting.
Specifically, the inlet conditions are assumed to be constant, reflecting the typical conditions observed in a direct-connect experimental setup.

% \subsubsection{Analytical Model}
% \todoUMBC{Should we delete this subsection since we are not including analytical model?}

% The control of the analytical model under various unknown operating conditions is described in detail in \cite{oveissi2023learning}, where it is shown that the adaptive controller is robust to operating conditions and uncertainties in the SFRJ model. 
% Furthermore, it is shown that the adaptive controller outperforms a fixed-gain controller. 

% \subsection{Open-loop Simulations}

\subsubsection{Operational Envelope}

\begin{figure}[h]
    \centering
    \includegraphics[width=3.25 in]{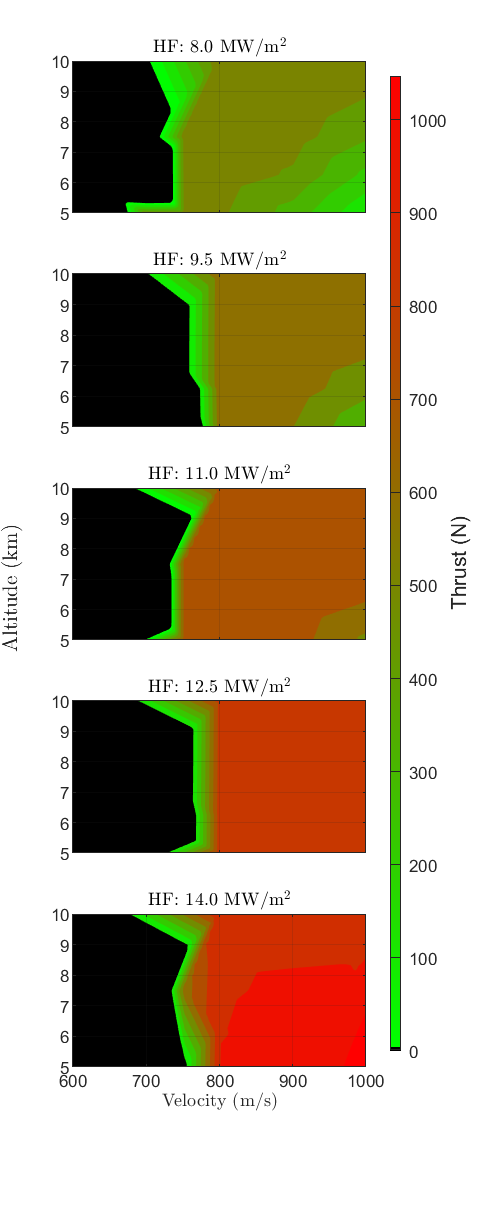}
    \caption{}
    \label{fig:1by5_OL}
\end{figure}

This section examines the influence of boundary conditions on SFRJ performance.
Specifically, the total inlet pressure, total inlet temperature, and heat flux are varied. 
The inlet pressure and temperature are determined using a standard atmospheric model for altitudes between $5 \ \rm km$ and $10\ \rm km$.
Figure \ref{fig:1by5_OL} shows the thrust generated by the SFRJ in various operating conditions. 
Note that the thrust values below zero are shown in black, and the heat flux, expressed in megawatts per square meter, is indicated in the title of each subplot.

% \todoUMBC{Alex - add OL results here}

% \begin{figure}
%     \centering
%     \begin{subfigure}[t]{0.3\textwidth}
%         \includegraphics[width=3.25 in]{Figures/HF1_OL_sweep.eps}    
%     \end{subfigure}
%     ~
%     \begin{subfigure}[t]{0.3\textwidth}
%         \includegraphics[width=3.25 in]{Figures/HF2_OL_sweep.eps}    
%     \end{subfigure}
%     ~
%     \begin{subfigure}[t]{0.3\textwidth}
%         \includegraphics[width=3.25 in]{Figures/HF3_OL_sweep.eps}    
%     \end{subfigure}
%     \\
%     \begin{subfigure}[t]{0.3\textwidth}
%         \includegraphics[width=3.25 in]{Figures/HF4_OL_sweep.eps}    
%     \end{subfigure}
%     ~
%     \begin{subfigure}[t]{0.3\textwidth}
%         \includegraphics[width=3.25 in]{Figures/HF5_OL_sweep.eps}    
%     \end{subfigure}
%     \caption{Caption}
%     \label{fig:enter-label}
% \end{figure}

% \begin{figure}[h]
%     \centering
%     \includegraphics[width=3.25 in]{Figures/HF1_OL_sweep.eps}
%     \caption{}
%     \label{fig:HF1_OL_sweep}
% \end{figure}

% \begin{figure}[h]
%     \centering
%     \includegraphics[width=0.6\columnwidth]{Figures/HF2_OL_sweep.eps}
%     \caption{}
%     \label{fig:HF2_OL_sweep}
% \end{figure}

% \begin{figure}[h]
%     \centering
%     \includegraphics[width=0.6\columnwidth]{Figures/HF3_OL_sweep.eps}
%     \caption{}
%     \label{fig:HF3_OL_sweep}
% \end{figure}

% \begin{figure}[h]
%     \centering
%     \includegraphics[width=0.6\columnwidth]{Figures/HF4_OL_sweep.eps}
%     \caption{}
%     \label{fig:HF4_OL_sweep}
% \end{figure}

% \begin{figure}[h]
%     \centering
%     \includegraphics[width=0.6\columnwidth]{Figures/HF5_OL_sweep.eps}
%     \caption{}
%     \label{fig:HF5_OL_sweep}
% \end{figure}

% \begin{figure}[h]
%     \centering
%     \includegraphics[width=3.25 in]{Figures/HF_total_OL_sweeps.eps}
%     \caption{}
%     \label{fig:HF_total_OL_sweeps}
% \end{figure}

\subsubsection{Thrust Regulation in Static Conditions.}

% First, we consider the problem of regulating the generated thrust in static conditions, that is, the boundary conditions at the inlet and exit are assumed to be constant. 
% This scenario is motivated by the ground test of engines. 

\begin{figure}[h]
    \centering
    \includegraphics[width=3.25 in]{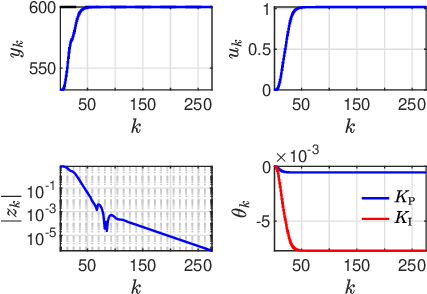}
    \caption{Closed-loop response of the SFRJ to a step command.}
    \label{fig:SU2_SFRJ_M3b_NominalHeat_10e6_NoThrustRef_steps_275_zfac_0_1_Nu_1_R0_1e6_linmap_1e6_ref_600_maxiter_1000}
\end{figure}

\textbf{Closed-loop Simulation and Hyperparameter Tuning.}
The SFRJ is commanded to generate a constant thrust value of $r = 600 $ N. 
The adaptive control law, consisting of adaptive gains $K_\rmp$ and $K_\rmi$, is optimized by the RCAC algorithm.
To tune the hyperparameters $N_1$ and $P_0$ of the RCAC algorithm, $N_1$ is fixed at $1$ and $P_0$ is varied logarithmically over the range $10^{-10}I_2$ to $10^{10}I_2.$
The value of $P_0$ is selected based on achieving visually acceptable closed-loop performance.
After preliminary tuning of the RCAC hyperparameters, we set $P_0 = 10^{-6}I_2, N_1 = 1$.
% \todoUMBC{Parham, are we using w0? Parham: yes}
Figure \ref{fig:SU2_SFRJ_M3b_NominalHeat_10e6_NoThrustRef_steps_275_zfac_0_1_Nu_1_R0_1e6_linmap_1e6_ref_600_maxiter_1000} shows the closed-loop response, where the first subplot shows the commanded and the generated thrust, the second subplot shows the control signal $u_k$ generated by RCAC, the third subplot shows the absolute value of the output error $z_k$ on a log scale, and the fourth subplot shows the PI controller gains $\theta_k$ updated by RCAC at each step.
Note that the third subplot illustrates the exponential convergence of the output error to zero and the exponential stability of the closed-loop system. 
We emphasize that the RCAC algorithm optimizes the controller coefficients using only the measured data and does not rely on the SFRJ model to optimize the controller. 

% \begin{table}[h]
%     \centering
%     \begin{tabular}{|l|c|}   
%         \hline
%         Case    & Figure \\ \hline
%         Constant step constant & Figure \ref{fig:SU2_SFRJ_M3b_NominalHeat_10e6_NoThrustRef_steps_275_zfac_0_1_Nu_1_R0_1e6_linmap_1e6_ref_600_maxiter_1000} \\ \hline
%         Sequence of step commands & Figure \ref{fig:SU2_SFRJ_M3b_NominalHeat_10e6_NoThrustRef_steps_600_zfac_0_1_Nu_1_R0_1e6_linmap_1e6_ref_600_700_800_maxiter_1000} \\ \hline
%         Effect of RCAC hyperparameters & Figure \ref{fig:SU2_SFRJ_M3b_Sensitivity_tests_v2} \\ \hline
%         Effect of Inlet Mach number, T, and P & vary M T P, keep thrust constant \\ \hline
%         % Effect of Inlet Temperature & \\ \hline
%         % Effect of Inlet Pressure & \\ \hline 
%         Comparison with fixed-gain controller & \\ \hline 
%     \end{tabular}
%     \caption{Caption}
%     \label{tab:my_label}
% \end{table}

% \textbf{idea 1 - run MC sample for M T and P}

% \todoUMBC{@Parham, is this figure updated with the latest map?  Parham: YES}

\begin{figure}[h]
    \centering
    \includegraphics[width=3.25 in]{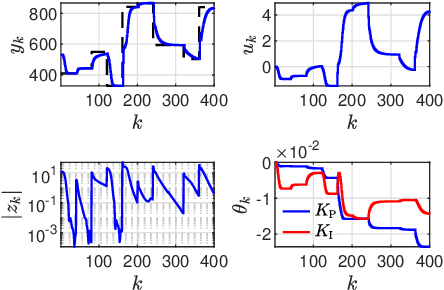}
    \caption{Closed-loop response of the SFRJ to a sequence of random step commands.}
    \label{fig:SU2_SFRJ_M3b_NominalHeat_10e6_NoThrustRef_steps_400_zfac_0_1_Nu_1_R0_1e6_linmap_1e6_ref_random_maxiter_1000}
\end{figure}

\textbf{Adaptation in Different Commands.}
Next, we investigate the performance of the adaptive controller when the reference commands vary over time. 
In this experiment, the reference signal is constructed as a sequence of randomly generated step commands within the range $[300, 900]$~N. The command is updated every 40 time steps to simulate varying operating conditions. Specifically, the reference at each interval is sampled from a uniform distribution:
\begin{align}
    r_k \sim \mathcal{U}(300, 900),
\end{align}
where $r_k$ is the constant reference value applied during the $k$-th interval.
It is important to emphasize that the RCAC hyperparameters are not re-tuned for this experiment.
% In RCAC, we set $P_0 = 10^{-6}I_2, N_1 = 1$.  
Figure~\ref{fig:SU2_SFRJ_M3b_NominalHeat_10e6_NoThrustRef_steps_400_zfac_0_1_Nu_1_R0_1e6_linmap_1e6_ref_random_maxiter_1000} shows the closed-loop response of the SFRJ system under this sequence of randomized reference inputs.
Same performance metrics as shown in Figure \ref{fig:SU2_SFRJ_M3b_NominalHeat_10e6_NoThrustRef_steps_275_zfac_0_1_Nu_1_R0_1e6_linmap_1e6_ref_600_maxiter_1000} are reported.
% 
% Note that the adaptive gains of the controller readjust as the command changes suggesting the adaptive algorithm is optimizing the controller accordingly. 
Note that the adaptive controller gains readjust in response to changes in the command signal, indicating that the adaptive algorithm is actively optimizing the controller to accommodate the evolving operating requirements.
% The first subplot shows the commanded and the generated thrust,
% the second subplot shows the control $u_k$ given by RCAC,
% the third subplot shows the absolute value of the output error $z_k$ on a log scale, 
% and the fourth subplot shows the PI controller gains $\theta_k$ updated by RCAC at each step.
% % 
% Note that the RCAC hyperparameters are not readjusted. 
% \begin{figure}[h]
%     \centering
%     \includegraphics[width=0.6\columnwidth]{Figures/SU2_SFRJ_modelV2_NominalHeat_2e6_steps_600_zfac_0_1_Nu_1_R0_1e8_ref_600_700_maxiter_1000.eps}
%     \caption{Closed-loop response of the SFRJ to a sequence of step commands}
%     \label{fig:SU2_SFRJ_modelV2_NominalHeat_2e6_steps_600_zfac_0_1_Nu_1_R0_1e8_ref_600_700_maxiter_1000}
% \end{figure}

% \begin{figure}[h]
%     \centering
%     \includegraphics[width=3.25 in]{Figures/SU2_SFRJ_M3b_NominalHeat_10e6_NoThrustRef_steps_600_zfac_0_1_Nu_1_R0_1e6_linmap_1e6_ref_600_700_800_maxiter_1000.eps}
%     \caption{Closed-loop response of the SFRJ to a sequence of step commands.}
%     \label{fig:SU2_SFRJ_M3b_NominalHeat_10e6_NoThrustRef_steps_600_zfac_0_1_Nu_1_R0_1e6_linmap_1e6_ref_600_700_800_maxiter_1000}
% \end{figure}

% \todoUMBC{@Parham, is this figure updated with the latest map? Parham: YES}

\begin{figure}[h]
    \centering
    \includegraphics[width=3.25 in]{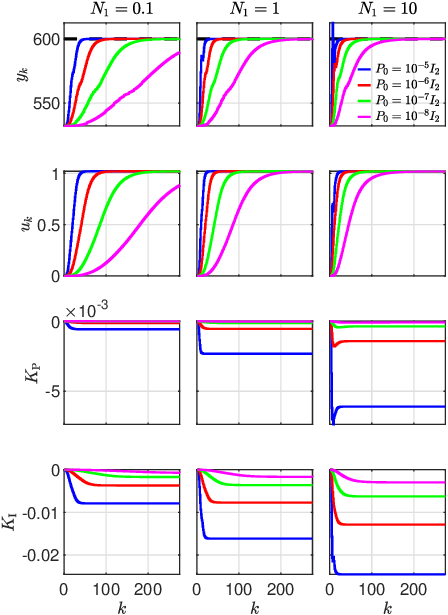}
    \caption{Effect of RCAC hyperparameters $P_0$ and $N_1$ on the closed-loop response of the SFRJ.}
    \label{fig:SU2_SFRJ_M3b_Sensitivity_tests_v3}
\end{figure}

\textbf{Robustness to Hyperparameter Choice.}
Next, the impact of the RCAC hyperparameters on closed-loop performance is investigated by reconsidering the step command-following problem.
% 
% We reconsider the step command-following problem.
% 
In RCAC, the hyperparameters are set as $N_1 = n$ and $P_0 = p I_2,$ where $n \in \{ 0.1, 1, 10\}$ and $p = \{10^{-5}, 10^{-6}, 10^{-7}, 10^{-8}\}.$
% Note that $I_2$ is the $2\times 2$ identity matrix. 
% The closed-loop simulation is thus run twelve times with all combinations of $n$ and $p.$
This results in a total of twelve closed-loop simulations, corresponding to all combinations of $n$ and $p.$
Figure \ref{fig:SU2_SFRJ_M3b_Sensitivity_tests_v3} shows the effect of RCAC hyperparameters on the closed-loop response of the SFRJ. 
The first row shows the thrust output $y_k$ of the SFRJ, and the second row shows the control $u_k$ used to generate the heat flux.
Finally the third and the fourth rows show the controller gains being updated. 
Note that a larger value of $P_0$ yields faster convergence but results in a larger overshoot. 
Similarly, a larger value of $N_1$ yields a faster response.

% Furthermore, to quantify the performance of the adaptive controller, we compare its performance with a fixed-gain controller, where the fixed gains are selected to be the gains optimized by RCAC at the end of the single step command problem shown in Figure \ref{fig:SU2_SFRJ_M3b_NominalHeat_10e6_NoThrustRef_steps_275_zfac_0_1_Nu_1_R0_1e6_linmap_1e6_ref_600_maxiter_1000}. 
 %

\textbf{Robustness to Operating Conditions.}
To evaluate the robustness of the adaptive control system, the effect of variations in inlet velocity, temperature, and pressure on the closed-loop thrust regulation performance is investigated. 
The inlet conditions are independently sampled from uniform distributions with the velocity in the range of 800–1000 $\rm m/s$, the temperature in the range of 220–260 $\rm K$, and the pressure in the range of 26–54 $\rm kPa$. 
A total of 15 samples are generated, each representing a distinct operating condition. 
For each case, the SFRJ is commanded to maintain a constant thrust of 600 N, regardless of the variations in inlet conditions.

\begin{figure}[h]
    \centering
    \begin{subfigure}[b]{0.49\textwidth}
        \centering
        \includegraphics[width=3.25 in]{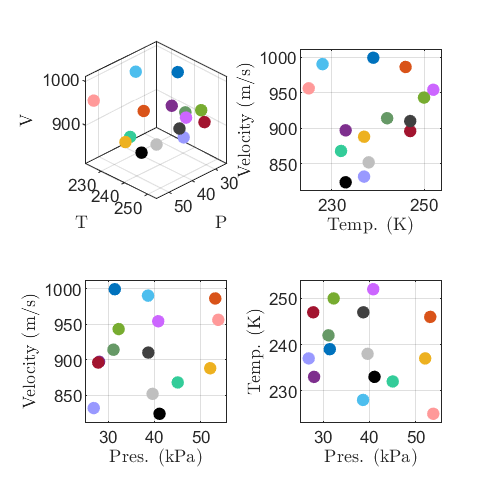}
        \caption{Distribution of the sampled inlet conditions.}
        \label{fig:SampleDistro}
    \end{subfigure}
    % \hspace{0.01\textwidth} % small space instead of \hfill
    ~
    \begin{subfigure}[b]{0.49\textwidth}
        \centering
        \includegraphics[width=3.25 in]{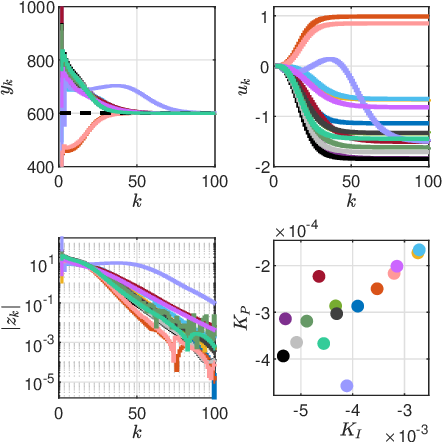}
        \caption{Closed-loop response.}
        \label{fig:M3bO1MC}
    \end{subfigure}
    \caption{Effect of inlet velocity, temperature, and pressure on the closed-loop performance of the SFRJ.}
    \label{fig:M3bO1MC_both}
\end{figure}

Figure \ref{fig:SampleDistro} shows the distribution of the sampled inlet conditions used in this study. 
Figure \ref{fig:M3bO1MC} shows the corresponding closed-loop thrust responses. Note the bottom right subplot shows the converged gains from the adaptive controller for each random inlet condition.
These results demonstrate that, despite significant uncertainty in the inlet conditions, the adaptive controller successfully regulates the thrust to the desired level, highlighting the robustness of the proposed control strategy.

\begin{figure}[h]
    \centering
    \begin{subfigure}[b]{0.49\textwidth}
        \centering
        \includegraphics[width=3.25 in]{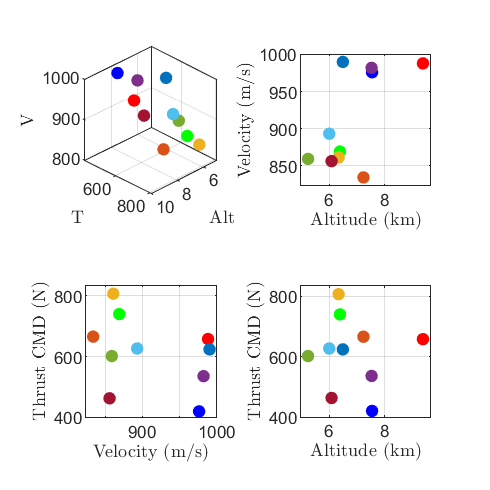}
F        \caption{Distribution of the sampled inlet conditions.}
        \label{fig:SampleDistro2}
    \end{subfigure}
    \hspace{0.01\textwidth} % small space instead of \hfill
    \begin{subfigure}[b]{0.49\textwidth}
        \centering
        \includegraphics[width = 3.25 in]{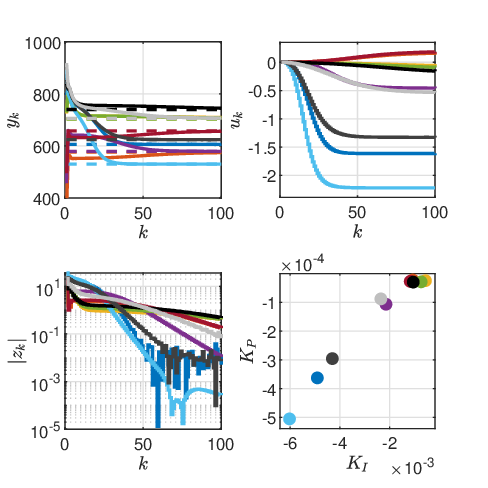}
        \caption{Closed-loop response.}
        \label{fig:M3bO1MC2}
    \end{subfigure}
    \caption{Effect of inlet velocity, altitude at various thrust commands on the closed-loop performance of the SFRJ.}
    \label{fig:M3bO1MC2_both}
\end{figure}

Figure \ref{fig:SampleDistro2} shows the distribution of sampled inlet conditions. Altitude and velocity were each drawn from a uniform distribution over a range of 5-10 km and 800-1000 $\rmm/\rms$, respectively. Thrust commands were then generated pseudo-randomly between 400-800 N, subject to the engine’s operational limits. Specifically, for each altitude–velocity pair we drew a candidate thrust value in 400-800 N and checked it against the precomputed input–output lookup table from Figure \ref{fig:1by5_OL}. If the candidate exceeded the maximum achievable thrust at that inlet condition, we resampled until a feasible thrust was obtained.

% \begin{figure}[h]
%     \centering
%     \includegraphics[width=3.25 in]{Figures/SampleDistro.eps}
%     \caption{Distribution of the sampled inlet conditions.}
%     \label{fig:SampleDistro}
% \end{figure}

% \begin{figure}[h]
%     \centering
%     \includegraphics[width=3.25 in]{Figures/M3bO1MC.eps}
%     \caption{Effect of inlet velocity, temperature, and pressure on the closed-loop performance of the SFRJ.}
%     \label{fig:M3bO1MC}
% \end{figure}

% \begin{figure}[h]
%     \centering
%     \includegraphics[width=3.25 in]{Figures/SampleDistro2.eps}
%     \caption{Distribution of the sampled inlet conditions.}
%     \label{fig:SampleDistro}
% \end{figure}

% \begin{figure}[h]
%     \centering
%     \includegraphics[width=3.25 in]{Figures/M3bO1MC2.eps}
%     \caption{Effect of inlet velocity, altitude at various thrust commands on the closed-loop performance of the SFRJ.}
%     \label{fig:M3bO1MC}
%     % The altitude is sampled from a uniform distribution between 5km and 10km, Thrust command sampled between uniformly distributed between 400 and 800 N, and velocity sampled between 800 mps and 1000 mps. 
% \end{figure}

\subsection{Dynamic Operating Conditions}
In this section, we consider a real-world deployment scenario in which the solid fuel ramjet (SFRJ) operates under dynamic flight conditions. 
Specifically, we assume the SFRJ serves as the primary propulsion system of a missile engaged in an interception maneuver.

As the missile is guided through the atmosphere by its onboard guidance and flight control systems, the SFRJ encounters time-varying boundary conditions. 
These variations lead to fluctuations in thrust output, even when the heat flux remains constant. 
Since most missile guidance algorithms assume constant thrust to compute the required normal acceleration for interception, the SFRJ control system must regulate the thrust to maintain a consistent output despite changing environmental conditions.

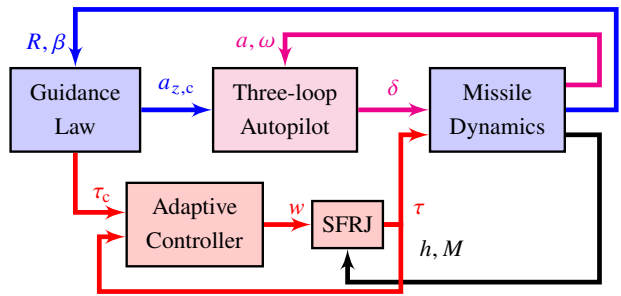
\begin{figure}[h]
    \centering
    \resizebox{3.25 in}{!}
    {
    \begin{tikzpicture}[auto, node distance=2cm,>=latex',text centered, line width = 1]

        \node [smallblock, fill=blue!20] at (-5,0) (Guidance) 
        {$\begin{array}{c}{\rm Guidance} \\ {\rm Law}\end{array}$};
        
        % \node at (-3,0) (reference) {$a_{z,\rm ref}$};
 
        \node [smallblock, fill=magenta!20, right = 3 em of Guidance] (Controller) 
        {$\begin{array}{c}{\text{Three-loop} } \\ {\rm Autopilot}\end{array}$};

        % \node [smallblock, fill=blue!20, right = 3 em of Controller] (Nonlinear) 
        % {$\begin{array}{c}{\rm Actuator} \\ {\rm Dynamics}\end{array}$};

        \node [smallblock, right = 3 em of Controller] (Plant) 
        {$\begin{array}{c}{\rm Missile} \\ {\rm Dynamics}\end{array}$};

        \node [smallblock, fill=red!20, below right = 2 em and -2 em of Controller] (SFRJ) {SFRJ};

        \node [smallblock, fill=red!20, left = 2 em of SFRJ] (RCAC) 
        {$\begin{array}{c}{\text{Adaptive} } \\ {\rm Controller}\end{array}$}
        % {Adaptive Controller}
        ;

        % \node[sum, left = 1.979 em of RCAC] (sum) {};
        
        % \node[left = 1.5 em of RCAC] (T_ref) {$T_{\rm ref}$};

        % \node [smallblock, right = 3 em of Plant] (IMU) {IMU};
        
        % \node[right = 3 em of IMU] (output) {$y$};

        \draw[->, blue, line width = 2] (Guidance) -- node[xshift = 0em, yshift = .2em]{$a_{z,\rm c}$} (Controller);
        
        % \draw[->] (reference) -- node[xshift = 0em, yshift = .2em]{} (Controller);
        
        \draw[->, magenta, line width = 2] (Controller)-- node[xshift = 0em, yshift = .2em]{$\delta$} (Plant);

        % \draw[<-] (Guidance.100)-- +(0,1) node[xshift = -1.25em, yshift = -1em]{$R, \beta$} ;

        \draw[->, black, line width = 2] (Plant.-20) -| +(0.5,-2.2) -|(SFRJ.270) node[xshift = 4em, yshift = -0em] {$h, M$};
        % \draw[->, thick] (Plant.232) |- +(0,-0.4) -| (RCAC.90);

        \draw[->, red, line width = 2] (SFRJ.east) node[xshift = 1.5em, above]{$\tau$} -| ++(0.25, -1) |- ++(-4.5,0) |- (RCAC.190);

        \draw[->, red, line width = 2] (SFRJ.east) -| ++(0.25,00) |- (Plant.200);

        % \draw[->] (SFRJ.west) -| (sum.north) node[near end, xshift = -1.3em , yshift = -1.3em]{$-$};
        
        \draw[->, red, line width = 2] (Guidance.south) |- (RCAC.170) node[xshift = -1em, above] {$\tau_{\rm c}$};

        % \draw[->] (sum.east) node[xshift = 1.0em, above]{$e$} |- (RCAC.west) ;

        \draw[->, red, line width = 2] (RCAC.east) node[xshift = 1.4em, above]{$w$} --(SFRJ) ;
       
        \draw[->, magenta, line width = 2] (Plant.20) -| ++(0.5,+0.75) -| (Controller.north) node[xshift = -1.25em, yshift = 1em]{$a, \omega$};
        \draw[->, blue, line width = 2] (Plant.-00) -| ++(0.75,+1.5) -| (Guidance.90) node[xshift = -1.25em, yshift = 1em]{$R, \beta$};
        
        % \draw[->] (IMU.east) node[xshift = 1em, yshift = 0.9em]{$y$} -- +(1,0) --  +(+1,-3.5)  -| (Controller.270);

    \end{tikzpicture}
    }
    % \vspace{-2em}
    \caption{Block diagram illustrating the integration of the guidance system, flight controller, engine controller, and missile dynamics. }
    \label{fig:integratedFramework}
\end{figure}

For this study, we assume that the missile engagement occurs within a vertical plane.
This simplification allows us to streamline the guidance law and flight control design without loss of generality. 
The equations governing the missile's longitudinal dynamics, along with the three-loop autopilot flight controller, are detailed in our previous work \cite{oveissi2023learning}.
Figure \ref{fig:integratedFramework} illustrates the integration of the guidance system, flight controller, engine controller, and missile dynamics.
The guidance law, based on the range $R$ and the line-of-sight angle $\beta,$ generates the normal acceleration command $a_{z, \rmc}$.
In this study, we assume that the guidance law is designed under the assumption of constant missile thrust.
The three-loop autopilot uses the normal acceleration measurement $a_z,$ the pitch rate measurement $\omega,$ and the command $a_{z, \rmc}$ to compute the required fin deflection angle $\delta$.
The constant thrust command is passed to the adaptive engine controller, which regulates the thrust produced by the SFRJ using the RCAC algorithm described earlier.

% \subsubsection{Inteception Dynamics}

% \subsubsection{Missile Dynamics}

% Missile dynamics and autopilot is described in \cite{oveissi2025learning}

% \subsubsection{Intercetpion Scenario}

% We consider an interception scenario where the missile is required to intercept an evader. 
% In all of the simulations considered in this paper, 

\begin{figure}[]
    \centering
    \includegraphics[width=0.5\linewidth]{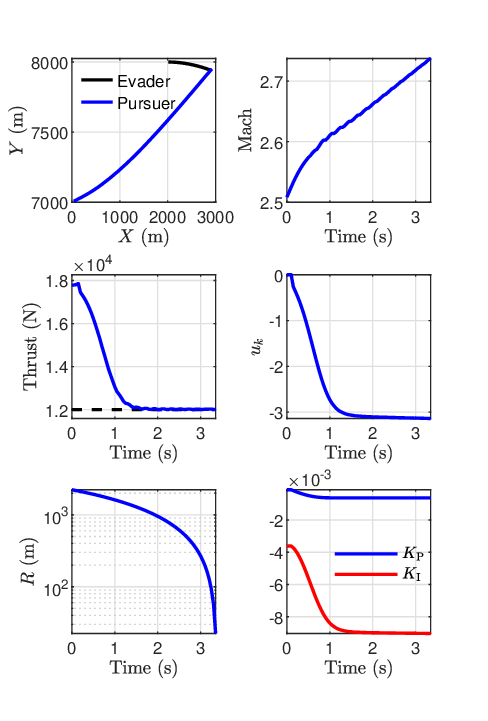}
    \caption{Interception trajectory of the missile and the closed-loop performance of the SFRJ.}
    \label{fig:M3b_M1_O3}
\end{figure}

In this study, the evader is modeled as having a mass of $10,000$ $\rm kg$ and a constant thrust of $76,310$ $\rm N.$
% and a normal acceleration capability of $10$ $\rm m/s^2.$
At the start of the simulation, the evader is assumed to be flying at Mach $0.75$ with a $0$ $\rm deg$ flight path angle at an altitude of $8$ $\rm km$ and at a horizontal distance of $2$ $\rm km$ from the pursuing missile. 
The pursuer is assumed to have a mass of $204$ $\rm kg$ and flying at Mach $2.5$ at a $10$ $\rm deg$ flight path angle with an initial angle of attack of $1$ $\rm deg$ at an altitude of $7$ $\rm km.$
The thrust command is set to $12$ $\rm kN.$
Figure \ref{fig:M3b_M1_O3} shows the interception trajectory of the missile and the closed-loop performance of the SFRJ.
The top row shows the trajectories of the evader and the pursuer. 
The first subplot in the middle row displays the commanded thrust and the thrust generated by the SFRJ, while the second subplot shows the adaptive control signal produced by RCAC.
In the bottom row, the first subplot shows the distance to the target on a logarithmic scale, and the second subplot shows the evolution of the controller gains as the engine experiences time-varying inlet conditions.

% In RCAC, we set $P_0$ to $1e-7$, $R_z$ to $1$, $R_u$ to $0$, and $N_1$ to $1$. 
% The thrust output of the SFRJ is scaled by 25 and 150 for M3b and M1 respectively to match the thrust requirements of the interceptor.  \todoUMBC{Alex, fill in the details here as before.}

% Consider the scenario where a pursuer has initial conditions of mach 1, 5 degree flight path angle, 5 degree pitch angle, and 1 km altitude. 

% The evader has an intial conditions of mach 1, 45 degree flight path angle, altitude of 5 km, and 5 km horizontal distance from the pursuer. The evader undergoes 10 $\frac{m}{s^2}$  of normal acceleration and a thrust of 76310 newtons with a mass of 10000 kg. 

%The SFRJ is ran in 3 cases. In case 1, the thrust command is set constant through the run time, that being 2 kN, in case 2, the thrust command is set to 2.1 kN for the first 25 seconds and 3 kN for anytime after, and in case 3, the thrust is commanded to be 2.2 kN for the first 16.6 seconds, 3.2 kN for the next 16.6 seconds, and 4 kN for anytime after that. 

% 
% \subsection{Unstart Control}

\iffalse
\section{Discussion}
In future work, more sophisticated combustion models are planned for implementation using the SU2-NEMO (NonEquilibrium MOdels) code to achieve a more accurate representation of the combustion process in the SFRJ \cite{maier2023}.
SU2-NEMO solves the Navier–Stokes equations for multi-species gases in thermochemical nonequilibrium, enabling the modeling of continuum hypersonic flows. The Mutation++ library \cite{scoggins2020a} (Multicomponent Thermodynamic And Transport properties for IONized gases in C++) supplies algorithms for computing thermodynamic properties and chemical kinetics of gas mixtures. 
SU2-NEMO is coupled with Mutation++ to enable detailed modeling of combustion chemistry relevant to SFRJs \cite{deboskey2023}. 
% These future developments are further discussed in Sec.~\ref{sec:conclusion}.
\fi

\section{Conclusions}
\label{sec:conclusions}

This paper presented a computational and control framework for thrust regulation in solid-fuel ramjets.
A CFD-based model with simplified combustion modeling was developed to characterize the thrust response and to establish the operational envelope, including the onset of inlet unstart.
Building upon this predictive model, an adaptive proportional–integral controller was designed and continually updated using the retrospective cost adaptive control algorithm.
The closed-loop simulations demonstrated that RCAC is capable of regulating thrust in both static and dynamic operating conditions without requiring an explicit model of the SFRJ dynamics.

The results highlight three key contributions. 
First, the CFD-based approach enabled the identification of thrust limits and unstart boundaries with sufficient fidelity to inform control design at manageable computational cost.
Second, the integration of RCAC with an adaptive PI controller provided effective thrust regulation under variations in reference commands, hyperparameters, and inlet conditions, thereby underscoring the robustness of the learning-based approach.
Third, dynamic simulations illustrated the ability of the adaptive controller to maintain commanded thrust during engagement scenarios with time-varying operating conditions, a critical requirement for practical deployment in air-breathing propulsion systems.

Overall, the study demonstrates that learning-based adaptive control, and RCAC in particular, offers a promising pathway for thrust regulation in SFRJs, where conventional model-based control strategies are hindered by strong nonlinearities, parametric uncertainty, and limited observability.
% Future work will extend this framework by incorporating higher-fidelity combustion models, investigating real-time implementation issues, and validating the approach through hardware-in-the-loop and experimental testing.
Importantly, the results suggest that adaptive, data-driven control architectures could play a central role in enabling reliable SFRJ operation for next-generation missile propulsion and long-range hypersonic flight systems.

\section{Acknowledgment}
This research was supported by the Office of Naval Research grant N00014-23-1-2468.
This research was supported in part through computational resources and services provided by the University of Arizona's Research Data Center (RDC). 
The authors would like to thank Brian Reitz and Alireza Farahmandi from NAWCWD China Lake for productive discussions on SFRJ physics.

% \section{Data-driven Learning Control}
% \section{Discussion}

% 
\bibliography{Refs}
\end{document}

%% file: PackagesCommands.tex
\usepackage{amsmath} % assumes amsmath package installed
\usepackage{amsfonts}
\usepackage{float} 
\usepackage{subcaption}

\usepackage{tikz}
\usetikzlibrary{shapes,arrows,calc,positioning}
\tikzstyle{bigblock} = [draw, fill=blue!20, rectangle, 
    minimum height=6em, minimum width=8em]
\tikzstyle{medblock} = [draw, fill=blue!20, rectangle, 
    minimum height=4em, minimum width=4em]    
\tikzstyle{mux} = [draw, fill=black!20, rectangle, 
    minimum height=5em, minimum width=0.1em]    
\tikzstyle{smallblock} = [draw, fill=blue!20, rectangle, 
    minimum height=2em, minimum width=3em]
    
\tikzstyle{data_block} = [draw, fill=green!20, rectangle, 
    minimum height=2em, minimum width=3em]
\tikzstyle{ops_block} = [draw, fill=blue!20, rectangle, 
    minimum height=2em, minimum width=3em]    
\tikzstyle{est_block} = [draw, fill=red!20, rectangle, 
    minimum height=2em, minimum width=3em]    
    
\tikzstyle{sum} = [draw, fill=blue!20, circle, node distance=1cm,minimum height=0.5cm]
\tikzstyle{signal} = [coordinate]
\tikzstyle{pinstyle} = [pin edge={to-,thin,black}]
\tikzstyle{block} = [draw, fill=blue!20, rectangle, 
    minimum height=3em, minimum width=9em]
\tikzstyle{blockS} = [draw, fill=blue!20, rectangle, 
    minimum height=3em, minimum width=4em]    
\tikzstyle{input} = [coordinate]
\tikzstyle{output} = [coordinate]
\usetikzlibrary{matrix}

\usetikzlibrary{positioning}
\usetikzlibrary{math}

\usepackage{hyperref}
\usepackage{xcolor}
\hypersetup{
    colorlinks,
    linkcolor={blue!100!black},
    citecolor={blue!50!black},
    urlcolor={blue!80!black}
}

\newcommand{\bc}{\begin{center}}
\newcommand{\ec}{\end{center}}
\newcommand{\benum}{\begin{enumerate}}
\newcommand{\eenum}{\end{enumerate}}

\newcommand{\matl}{\left[ \begin{array}}
\newcommand{\matr}{\end{array} \right]}

\renewcommand{\matl}{\begin{bmatrix}}
\renewcommand{\matr}{\end{bmatrix}}

\newcommand{\matls}{\left[ \begin{smallmatrix}}
\newcommand{\matrs}{\end{smallmatrix} \right]}
\newcommand{\isdef}{\stackrel{\triangle}{=}}

\newcommand{\rmI}{{\rm I}}

\newcommand{\rmP}{{\rm P}}

\newcommand{\rmT}{{\rm T}}

\newcommand{\rmc}{{\rm c}}
\newcommand{\rmd}{{\rm d}}

\newcommand{\rmi}{{\rm i}}

\newcommand{\rmm}{{\rm m}}

\newcommand{\rmp}{{\rm p}}

\newcommand{\rms}{{\rm s}}
\newcommand{\rmt}{{\rm t}}

\newcommand{\rmv}{{\rm v}}

\newcommand{\SO}{{\mathcal O}}

\usepackage{enumitem,amssymb}
\newlist{todolist}{itemize}{2}
\setlist[todolist]{label=$\square$}
\usepackage{pifont}